\documentclass[12pt]{amsart}
\usepackage{a4wide}
\usepackage{amsfonts}
\usepackage{bbm}
\usepackage{hyperref}
\newcommand{\E}{\mathbb E}
\newcommand{\reals}{\mathbb R}
\newcommand{\ints}{\mathbb N}
\renewcommand{\P}{\mathbb P}

\newcommand{\1}{\mathbbm 1}

\renewcommand{\epsilon}{\varepsilon}
\newcommand{\half}{\mbox{$\frac 1 2$}}
\newcommand{\elem}{\mathcal E}
\newcommand{\simple}{\mathcal S}

\usepackage{color}

\title{Stochastic Integration in ordered Banach Spaces}
\author{Joris Bierkens}
\author{Onno van Gaans}
\date{\today}
\begin{document}
\maketitle
\section{Vector-valued random variables}
Throughout this section, let $X$ be a Banach space with dual $X^*$ and Borel $\sigma$-algebra $\mathcal X$, and let $(\Omega, \mathcal F, \P)$ be a probability space. We want to treat functions $x: \Omega \rightarrow X$ as random variables. To make this precise, we give the following definitions.

\subsection{Vector-valued random variables}
A \emph{vector valued random variable} is an $\mathcal F$-measurable function $x: \Omega \rightarrow X$, that is
\[ x^{\leftarrow} (A) \in \mathcal F \quad \textrm{for all} \ A \in \mathcal X. \]

The \emph{$\sigma$-algebra generated by maps $x_\gamma:\Omega \rightarrow X$, $\gamma \in C$} is $\sigma(x_\gamma^{\leftarrow}(A) : A \in \mathcal X, \gamma \in C )$ and denoted by $\sigma(x_\gamma: \gamma \in C)$.

Random variables (real- or vector-valued) $x_1, x_2, \hdots$ are \emph{independent} if the $\sigma$-algebras $\sigma(x_1), \sigma(x_2), \hdots$ are independent. 
\subsection{Expectation}
Let $x \in L^1(\Omega;X)$, i.e. $x$ is strongly integrable.\footnote{Refer to Appendix~\ref{app:integration} for details.} Define the \emph{(vector-valued) expectation of $x$} (notation $\E_X [x]$) as
\[ \E_X [x] := \int_\Omega x \ d \P.\]

Proposition~\ref{prop:bochner=pointwise} shows that this expectation works `as expected' for function spaces:
\[ \E_X[x(\cdot)] = \E[x](\cdot), \]
i.e. evaluating the vector-valued expectation amounts to a pointwise evaluation of the usual real-valued expectation. Therefore we will drop the subscript $X$ in $\E_X$ and use the symbol $\E$ for both vector-valued and real-valued expectation, as long as this does not lead to confusion.

We will now verify some properties of vector-valued expectation.
\subsection{Lemma}
\label{sec:indeplemma}
Let $x \in L^1(\Omega; X)$ and $r \in L^1(\Omega)$ be independent vector- and real-valued integrable random variables, respectively. Then

\begin{itemize}
\item[(i)]
The random variables $r$ and $\phi(x)$ are independent for all $\phi \in X^*$.
\item[(ii)] $ \E_X [ r x ] = \E[r] \E_X [x].$
\end{itemize}
\begin{quote}
\textsc{Proof:}
\begin{itemize}
\item[(i)]
Let $\mathcal B$ be the Borel $\sigma$-algebra of $\reals$. For $A \in \mathcal B$,
\[ (\phi(x))^{\leftarrow}(A) = x^{\leftarrow} (\phi^{\leftarrow} (A)) \in x^\leftarrow (\mathcal X),\]
since $\phi^\leftarrow (A) \in \mathcal X$ because $\phi$ is continuous and hence Borel measurable. Therefore $(\phi(x))^\leftarrow(\mathcal B) \subset x^\leftarrow(\mathcal X)$ and therefore $(\phi(x))^\leftarrow(\mathcal B) $ and $r^\leftarrow(\mathcal X)$ are independent, so $\phi(x)$ and $r$ are independent.
\item[(ii)]
For all $\phi \in X^*$,
\begin{align*} \phi ( \E_X[ r x ] ) &= \E [ \phi ( r x) ] = \E [ r \phi(x)] = \E [r] \E [\phi(x)] \\
& = \E[ r] \phi (\E_X[ x]) = \phi(\E[ r] \E_X[ x]),\end{align*}
using that $\phi(x)$ and $r$ are independent, and that
\[ \E \left[ |\phi(x) | \right] \leq \E \left[ ||\phi ||_{X^*}  || x ||_X \right] \leq ||\phi||_{X^*} \E \left[ ||x ||_X \right] < \infty ,\]
so that $\phi(x) \in L^1(\Omega)$.
\hfill $\square$
\end{itemize}
\end{quote}

In the remainder of this section, unless stated otherwise, $\mathcal F$ will denote a $\sigma$-algebra on $\Omega$ and $\mathcal G$ a sub-$\sigma$-algebra of $\mathcal F$.

\subsection{Theorem and Definition (Conditional Expectation)}
\label{thm:conditionalexpectation} 
Let $(\Omega, \mathcal F, \P)$ a probability space and $x \in L^1(\Omega;X)$. Then there exists a $y \in L^1(\Omega; X)$ such that
\begin{enumerate}
\item
$y$ is $\mathcal G$ measurable,
\item
For every $A \in \mathcal G$ we have
\[ \E [ y \1_A ] = \E [ x \1_A ].\]
\end{enumerate}
Moreover, if $\tilde y$ is another $X$-valued random variable with these properties then $\tilde y = y$, a.s..

The random variable $y$ with the above properties is called \emph{(a version of) the conditional expectation of $x$ given $\mathcal G$}, and written as
\[ y = \E [ x | \mathcal G ], \quad \mathrm{a.s.} \]

Furthermore, conditional expectation satisfies the following properties:
\begin{enumerate}
\item[(i)] If $\mathcal F = \mathcal G$ then $y = x$, a.s.
\item[(ii)] \emph{(Linearity)} For $x_1, x_2 \in L_1((\Omega, \mathcal F), X)$, if $y_1 = \E [x_1 | \mathcal G]$ and $y_2 = \E [ x_2 | \mathcal G]$, and $a, b \in \reals$, then $ \E [ a x_1 + b x_2 | \mathcal G ] = a y_1 + b y_2, \quad \textrm{a.s.} $
\item[(iii)] If $x \in X$ and $r: \Omega \rightarrow \reals$, then $ \E [ x r | \mathcal G] = x \E[r|\mathcal G]$. 
\end{enumerate}

\begin{quote}
\textsc{Proof:}

\emph{Uniqueness.} Let $x \in X$. Suppose both $y$ and $\tilde y$ satisfy the requirements. Then $\E[(y-\tilde y) \1_A ] = 0$ for all $A \in \mathcal G$. But if $\E[ ||y - \tilde y||_X ] \neq 0$, then there exists a set $A \in \mathcal G$ with $y - \tilde y \neq 0$ on $A$, which gives a contradiction.

\emph{Existence.}
First suppose that $x$ is simple, and of the form
\[ x = \sum_{i=1}^n x_i \1_{A_i}, \quad x_i \in X, A_i \in \mathcal F \ \textrm{disjoint}, i = 1, \hdots, n .\]
Define
\[ \E \left[ x | \mathcal G \right] := \sum_{i=1}^n x_i \E \left[ \1_{A_i} | \mathcal G \right].\]
It is easily verified that this construction satisfies the requirements.

If $(x^n)_{n \in \ints}$ is a sequence of simple functions that approaches $x$ in $L^1((\Omega, \mathcal F),X)$, then
\begin{align*}
\E \left[ \left|\left| \E \left[ x^n - x^m | \mathcal G  \right] \right|\right|_X \right] & = \E \left[ || \sum_i (x_i^n - x_i^m ) \E \left[ \1_{A_i} | \mathcal G \right] ||_X \right] \\
& \leq \sum_i \E \left[ \E [ \1_{A_i} | \mathcal G ] \right] ||x_i^n - x_i^m ||_X \\
& = \sum_i \P(A_i) ||x_i^n - x_i^m || = ||x^n - x^m ||_{L^1(\Omega,X)}.
\end{align*}
So $( \E \left[ x^n | \mathcal G \right])_{n \in \ints}$ is a Cauchy sequence in $L^1((\Omega, \mathcal G), X)$ and hence has a unique limit, that we denote by $\E\left[ x | \mathcal G \right]$.
\end{quote}

\subsection{Lemma (Tower property)}
\label{lem:towerproperty}
Let $\mathcal G \subset \mathcal H \subset \mathcal F$ be $\sigma$-algebras. Let $x \in L^1((\Omega,\mathcal F);X)$, then
\[ \E[ x | \mathcal G] = \E [ \E[x|\mathcal H] | \mathcal G], \quad \textrm{a.s.}.\]
\begin{quote}
\textsc{Proof:}
For $A \in \mathcal G$,
\[ \E [x \1_A] = \E [ \E [x|\mathcal H] \1_A] = \E [ \E[\E [x|\mathcal H] |\mathcal G] \1_A]. \]
\hfill $\square$
\end{quote}

\subsection{Lemma}
\label{lem:normestimatecondexp}
If $x \in L^1((\Omega, \mathcal F); X)$ and $\mathcal G$ is a sub-$\sigma$-algebra of $\mathcal F$, then
\[ \left|\left| \E \left[ x | \mathcal G \right] \right|\right|_X \leq \E \left[ ||x||_X \ | \mathcal G \right], \quad \textrm{a.s.}\]
\begin{quote}
\textsc{Proof:}
If $x$ is a simple function, then this is easily seen from the construction of the conditional expectation. If $x \in L^1((\Omega,\mathcal F), X)$ and $(x_n)_{n \in \ints}$ approximates $x$ in $L^1((\Omega,\mathcal F), X)$, then as in the proof of Theorem~\ref{thm:conditionalexpectation}, $\E \left[ x_n | \mathcal G \right] \rightarrow \E \left[ x | \mathcal G \right]$ as $n \rightarrow \infty$ in $L^1(\Omega;X)$, and 
\begin{align*} \E \left[ \left| \E \left[ ||x||_X \ | \mathcal G \right] - \E \left[ ||x_n ||_X \ | \mathcal G \right] \right| \right] & \leq \E \left[ \E \left[ \left| ||x||_X - ||x_n||_X \right| \ | \mathcal G \right] \right] \\
& \leq \E \left[ ||x-x_n||_X \right] \rightarrow 0, \quad n \rightarrow \infty.
\end{align*}
So $\E \left[ x_n | \mathcal G \right] \rightarrow \E \left[ x | \mathcal G \right]$ in probability and $\left|\left| \E \left[ x_n | \mathcal G \right] \right|\right|_X \rightarrow \left|\left| \E \left[ x | \mathcal G \right] \right|\right| _X$ in probability. But then, for any $\epsilon > 0$, it is impossible that
\[ \P \left( \left|\left| \E [ x | \mathcal G ] \right|\right|_X - \E \left[ ||x||_X | \mathcal G \right] \geq \epsilon \right) > 0,\]
and therefore
\[ \P \left( \left|\left| \E [ x | \mathcal G ] \right|\right|_X \leq \E \left[ ||x||_X | \mathcal G \right] \right) = 1.\]
\hfill $\square$
\end{quote}

\subsection{Lemma (Dominated convergence for conditional expectation)}
\label{lem:condexpalmostsureconvergence}
If $x \in L^1((\Omega, \mathcal F);X)$ and $(x_i)_{i \in \ints}: \Omega \rightarrow X$ are strongly $\mathcal F$-measurable functions such that $||x_i - x||_X \rightarrow 0$ almost surely as $i \rightarrow \infty$, and if there is a function $y \in L^1((\Omega;\mathcal F), \reals)$, $y\geq 0$ such that $||x_i||_X \leq y$ almost surely for all $i \in \ints$, then $x_i \in L^1((\Omega;\mathcal F);X)$ for all $i \in \ints$ and
\[ \left|\left| \E \left[ x_i |\mathcal G \right] - \E \left[ x |\mathcal G \right] \right|\right|_X \rightarrow 0 \quad \textrm{a.s. as} \ i \rightarrow \infty.\]
\begin{quote}
\textsc{Proof:}
The $x_i$ are in $L^1((\Omega;\mathcal F);X)$ because
\[ \E \left[ ||x_i|| \right] \leq \E \left[ y \right] < \infty.\]

The claim follows because, by Lemma~\ref{lem:normestimatecondexp}
\begin{align*} \left|\left| \E \left[ x_i |\mathcal G \right] - \E \left[ x |\mathcal G \right] \right|\right|_X & \leq \E \left[ ||x_i -x||_X |\mathcal G \right],
\end{align*}
and the latter expression converges to zero by dominated convergence for real-valued conditional expectation.
\hfill $\square$
\end{quote}

\subsection{Lemma (Take out what is known I)}
\label{lem:takeoutwhatisknowni}
Let $x \in L^1((\Omega,\mathcal G);X)$ and $r \in L^1((\Omega, \mathcal F);\reals)$ such that $x r \in L^1((\Omega, \mathcal F);X)$. Then
\[ \E \left[ x r | \mathcal G \right] = x \E \left[ r |\mathcal G \right]. \]
\begin{quote}
\textsc{Proof:}
For $x \in L^1((\Omega,\mathcal G);X)$ a simple function, the result is immediate.
Now let $x_n$ be a sequence of $\mathcal G$-measurable simple functions so that $x_n \rightarrow x$ in $L^1(\Omega;X)$. 

\textsc{Claim:}
\[ x_n r \1_{\left\{||x_n||_X \leq ||x||_X + 1\right\}} \rightarrow x r, \quad \textrm{a.s.}.\]

\begin{quote}
\textsc{Proof of claim:}
Let $\Omega_0 = \{ \omega: x_n(\omega) \rightarrow x(\omega) \ \textrm{as} \ n \rightarrow \infty\}$, so $\P(\Omega_0) = 1$.
Let $\omega \in \Omega_0$ and $\epsilon > 0$. Pick $N \in \ints$ so that for $n \geq N$, we have that $||x_n(\omega) -x(\omega)||_X < 1 \wedge \epsilon$.
So for $n \geq N$, 
\[ || x_n(\omega) \1_{\left\{||x_n||_X \leq ||x||_X + 1 \right\}}(\omega) - x(\omega) ||_X < \epsilon.\]
So $x_n \1_{||x_n||_X \leq ||x||_X + 1} \rightarrow x$ almost surely, and therefore $x_n r \1_{||x_n||_X \leq ||x||_X + 1} \rightarrow x r$ almost surely, and the claim is proven.
\end{quote}

We have seen that since $x_n$ is simple,
\[ \E [x_n \1_{\left\{||x_n||_X \leq ||x||_X + 1 \right\}} r |\mathcal G] = x_n \E [\1_{\left\{ ||x_n||_X \leq ||x||_X + 1\right\}} r |\mathcal G].\]
By Lemma~\ref{lem:condexpalmostsureconvergence}, since $||x_n \1_{\left\{||x_n||_X \leq ||x||_X + 1\right\}} r ||_X \leq |r| \left(||x|| + 1\right)$,
\[ \E [ x_n \1_{\left\{ ||x_n||_X \leq ||x||_X + 1 \right\}} r |\mathcal G] \rightarrow \E [ x r |\mathcal G].\]
Also 
\begin{align*} & \left|\left| x_n \E \left[ \1_{\left\{||x_n||_X \leq ||x||_X + 1\right\}} r |\mathcal G \right] - x \E \left[ r |\mathcal G \right] \right|\right|_X \\
& \leq \left|\left| x_n \E \left[ \1_{\{\hdots\}} r |\mathcal G \right] - x_n \E \left[ r |\mathcal G \right] \right|\right|_X + \left|\left| x_n \E \left[ r |\mathcal G \right] - x \E \left[ r |\mathcal G \right] \right|\right|_X \\
& = ||x_n||_X \left| \E \left[ (\1_{\left\{||x_n||_X \leq ||x||_X + 1\right\}} - 1) r |\mathcal G\right] \right| + ||x_n - x||_X \left|\E \left[ r |\mathcal G \right] \right|,
\end{align*}
of which both the first term (by real-valued dominated convergence of conditional expectation) and the second term tend to zero as $n \rightarrow \infty$ on $\Omega_0$.
We conclude that
\begin{align*}
\E \left[ x r | \mathcal G \right] & = \lim_{n \rightarrow \infty} \E \left[ x_n \1_{\left\{||x_n||_X \leq ||x||_X + 1\right\}} r |\mathcal G \right] \\
& = \lim_{n \rightarrow \infty} x_n  \E \left[ \1_{\left\{||x_n||_X \leq ||x||_X + 1\right\}} r |\mathcal G \right] \\
& = x \E \left[ r |\mathcal G \right] \quad \textrm{almost surely}.
\end{align*}
\hfill $\square$
\end{quote}

\subsection{Lemma (Independence and conditional expectation)}
\label{lem:condexpectationindependence}
Suppose $\mathcal F, \mathcal G$ are independent $\sigma$-algebras, and let $x \in L^1((\Omega,\mathcal F);X)$. Then
\[ \E \left[ x | \mathcal G \right] = \E [x] \quad \textrm{almost surely}.\]
\begin{quote}
\textsc{Proof:}
It is straightforward to show this in case $x$ is a simple function. In case $x \in L^1((\Omega,\mathcal F);X)$, let $(x_i)$ be a sequence of simple functions such that $x_i \rightarrow x$. Then $x_i \1_{\{||x_i||_{L^1(\Omega;X)} \leq ||x||_{L^1(\Omega;X)}+1 \} } \rightarrow x$ almost surely and by dominated convergence for conditional expectation the result follows.
\hfill $\square$
\end{quote}

\subsection{Lemma}
If $x \in L^1((\Omega, \mathcal F); X)$, $\mathcal G$ is a sub-$\sigma$-algebra of $\mathcal F$ and $\phi \in X^*$ then
\[ \phi \left( \E\left[ x | \mathcal G \right] \right) = \E \left[ \phi(x) | \mathcal G \right], \quad \textrm{a.s.}.\] 

\begin{quote}
\textsc{Proof:} Because $\E\left[ x | \mathcal G \right]$ is $\mathcal G$-measurable, so is $\phi \left( \E \left[ x | \mathcal G \right] \right)$. If $A \in \mathcal G$, then almost surely
\begin{align*} \E \left[ \E \left[ \phi (x) |\mathcal G \right] \1_A \right] &= \E\left[ \phi(x) \1_A \right] = \E \left[ \phi(x \1_A) \right] = \phi \left (\E \left[ x \1_A \right] \right) \\
& = \phi \left( \E \left[ \E \left[ x | \mathcal G \right] \1_A \right] \right) = \E \left[ \phi \left( \E \left[ x | \mathcal G \right]\right) \1_A \right],
\end{align*}
which completes the proof.
\hfill $\square$
\end{quote}

\subsection{Theorem (Conditional Jensen) \cite{ting}}
\label{thm:cjensen}
Let $X$ be a Banach space and $Y$ an ordered Banach space such that the positive cone is closed (e.g. a Banach lattice). Let $x \in L^1((\Omega, \mathcal F);X)$, let $\mathcal G$ be a sub-$\sigma$-algebra of $\mathcal F$ and let $f : \Omega \times X \rightarrow Y$ be a function, such that
\begin{itemize}
\item[(i)] $f(\cdot, z)$ is strongly $\mathcal G$-measurable for every $z \in X$;
\item[(ii)] $f(\cdot, x(\cdot)) \in L^1((\Omega,\mathcal F);Y)$;
\item[(iii)] For every $\omega \in \Omega$, the function $f(\omega,\cdot) : X \rightarrow Y$ is continuous, and
\begin{equation} \label{eq:convex} f(\omega, \lambda x_1 + (1-\lambda) x_2) \preceq \lambda f(\omega, x_1) + (1-\lambda) f(\omega, x_2), \quad x_1, x_2 \in X, \lambda \in [0,1],\end{equation}
i.e. $f(\omega,\cdot)$ is convex.
\end{itemize} 

Then
\[ f\left(\omega, \E[ x | \mathcal G](\omega) \right) \preceq \E \left[ f(\cdot, x(\cdot)) | \mathcal G\right](\omega), \quad \textrm{for almost all} \ \omega \in \Omega. \]

\begin{quote}
\textsc{Proof:} See \cite{ting}).
\end{quote}

\subsection{Corollary (Jensen)}
\label{cor:jensen}
Let $X$ and $Y$ be as in the previous theorem.
Let $(S, \Sigma, \mu)$ be a finite measure space. Let $x \in L^1((S,\mu); X)$ and let $f: X \rightarrow Y$ be a continuous, convex function such that $f(\mu(S) x) \in L^1((S,\mu);Y)$. Then 
\[ f\left( \int_S x \ d \mu  \right) \preceq \frac 1 {\mu(S)} \int_S f(\mu(S) x) \ d \mu.\]

\begin{quote}
\textsc{Proof:}
Define a probability measure on $S$ by letting $\P(A) := \mu(A) / \mu(S)$ for $A \in \Sigma$. Let $\tilde x := \mu(S) x$ and let $\mathcal G$ be the trivial $\sigma$-algebra on $S$. Theorem~\ref{thm:cjensen} states that
\[ f\left( \frac 1 {\mu(S)} \int_S \tilde x \ d \mu  \right) \preceq \frac 1 {\mu(S)} \int_S f(\tilde x) \ d \mu.\]
\hfill $\square$
\end{quote}

\section{The It\^o integral for Banach space valued stochastic processes}
\label{sec:stochasticintegral}

In this section we describe how to construct the It\^o integral for general Banach spaces, with the least possible assumptions. In the succeeding sections this theory will be applied to the Banach spaces of continuous functions and of $L^p$-functions.

\subsection{Definition (Multiplication)}
\label{def:multiplication}
Let $X$ be a Banach space and $Y$ be a Banach lattice with partial ordering $\preceq$.\footnote{See Appendix~\ref{app:Banachlattice} and below for details.}
A symmetric, bilinear mapping $\cdot: X \times X \rightarrow Y$ is called \emph{$Y$-valued multiplication of elements in $X$} if it satisfies, for all $x, y \in X$,
\begin{itemize}
\item[(i)] $ || x \cdot y ||_Y \leq ||x ||_X \ ||y||_X $;.
\item[(ii)] $x^2 := x \cdot x \succeq 0$;
\item[(iii)] $x^2 = 0$ if and only if $x = 0$.
\end{itemize}

\subsection{Examples}
\label{multiplicationexamples}
Examples (i) and (ii) involve function spaces, where a partial ordering is given by $f \preceq g$ if $f(x) \leq g(x)$ for (almost) all $x$ in the domain of $f$ and $g$.

All examples below satisfy the conditions of Definition~\ref{def:multiplication}.

\begin{itemize}
\item[(i)] Let $X = C(K)$, with $K$ a compact metric space, and define $\cdot: X \times X \rightarrow X$ to be pointwise multiplication. Then
\begin{align*}
|| x \cdot y ||_X = \sup_{k \in K} |x(k) y(k)| \leq \sup_{k \in K} |x(k)| \ \sup_{k \in K} |y(k)| = ||x||_X \ ||y||_X,
\end{align*}
and
\[ || x^2||_X = \sup_{k \in K} |x^2(k)| = \left( \sup_{k \in K} |x(k)| \right)^2 = ||x||_X^2.\]
\item[(ii-a)] Let $(S, \Sigma, \mu)$ be a $\sigma$-finite measure space, $X = L^p(\mu)$, $2 \leq p < \infty$. Define $\cdot: X \times X \rightarrow Y$ to be pointwise multiplication, with $Y = L^{\frac p 2}(\mu)$.
Then, by Cauchy-Schwarz,
\begin{align*} || x \cdot y||_Y & = \left(\int_S |xy|^{\frac p 2} \ d \mu \right)^{\frac 2 p} \\
& \leq \left(\int_S |x|^p \ d \mu \right)^{\frac 1 p} \left(\int_S |y|^p \ d \mu \right)^{\frac 1 p} = ||x||_X ||y||_X.
\end{align*}
In this case
\[ || x^2 ||_Y = \left( \int_S |x^2|^{p/2} \ d \mu \right)^{2/p} = || x ||_X^2.\]
\item[(ii-b)] Let $(S, \Sigma, \mu)$ be a $\sigma$-finite measure space, $X = Y = L^\infty(\mu)$, with $\cdot: X \times X \rightarrow X$ pointwise multiplication. As for $C(K)$,
\[ ||x \cdot y||_X \leq ||x||_X \ ||y||_X, \quad \textrm{and} \quad ||x^2||_X = ||x||_X^2.\]
\item[(iii)] Let $X$ be a separable Banach space. By Lemma~\ref{lem:normingsequence}, there exists a sequence $(\phi_n)_{n \in \ints}$ in $X^*$ with $||\phi_n||_{X^*} = 1$ for all $n \in \ints$, that is norming for $X$, i.e. 
\[ ||x||_X = \sup_{n \in \ints} |\phi_n(x)| \quad \textrm{for all} \ x \in X.\]
Define a multiplication $\cdot$ by
\[ \cdot : X \times X \rightarrow l^\infty, \quad (x \cdot y)_n = \phi_n(x) \phi_n(y).\]
Then 
\[ || x \cdot y ||_{l^\infty} = \sup_{n \in \ints} |\phi_n(x) \phi_n(y)| \leq \sup_{n \in \ints} |\phi_n(x)| \ \sup_{m \in \ints} |\phi_m(y)| = ||x||_X \ ||y||_X.\]
A partial ordering on $l^\infty$ is given by $a \preceq b$ if $a_n \leq b_n$ for all $n \in \ints$. Then
\[ x^2 = (\phi_n(x)^2)_{n \in \ints} \succeq 0.\]
Finally, if $x^2 = 0$, then 
\[ ||x||_X^2 = \sup_{n \in \ints} |\phi_n(x)|^2 = 0,\] so that $x = 0$.
\item[(iv)] Let $X$ be a Hilbert space with inner product $(\cdot,\cdot)$. Define a multiplication by
\[ \cdot : X \times X \rightarrow \reals, \quad x \cdot y = (x,y).\]
It is readily verified that this multiplication structure satisfies the conditions of Definition~\ref{def:multiplication}.
\item[(v)] If $X = Y$ is a Banach algebra and a Banach lattice such that
\begin{itemize}
\item[(a)] $x^2 \succeq 0$ for all $x \in X$, and
\item[(b)] $x^2 = 0$ if and only if $x = 0$ in $X$,
\end{itemize}
then the multiplication structure of $X$ satisfies the conditions of~\ref{def:multiplication} (see Example (i)).
\end{itemize}

In the remainder of this section, the pair of Banach spaces $(X,Y)$ will satisfy the conditions of Definition~\ref{def:multiplication}.

\subsection{Lemma (Take out what is known II)} 
\label{lem:takeoutwhatisknownii}
If $x \in L^1((\Omega, \mathcal G);X)$ and $y \in L^1((\Omega,\mathcal F);X)$ such that
$x \cdot y \in L^1((\Omega, \mathcal F);Y)$
then \[ \E \left[ x \cdot y | \mathcal G \right] = x \cdot \E \left[ y | \mathcal G \right].\]
\begin{quote}
\textsc{Proof:}
If $y$ is simple, say $y = \sum_{i=1}^n y_i \1_{A_i}$, with $A_i$ disjoint elements of $\mathcal F$ and $y_i \in X$ for $i=1,\hdots$, then
\begin{align*} \E \left[ x \cdot y | \mathcal G \right] & = \sum_{i=1}^n \E \left[ (x \cdot y_i) \1_{A_i} | \mathcal G \right] = \sum_{i=1}^n (x \cdot y_i) \E \left[ \1_{A_i} | \mathcal G \right] \\
& = x \cdot \left( \sum_{i=1}^n y_i \E\left[\1_{A_i} | \mathcal G \right] \right) = x \cdot \sum_{i=1}^n \E \left[ y_i \1_{A_i} | \mathcal G \right] \\
& = x \cdot \E \left[ y | \mathcal G \right]
\end{align*}
by applying Lemma~\ref{lem:takeoutwhatisknowni} twice.

Now let $y \in L^1(\Omega;X)$ such that $x\cdot y \in L^1(\Omega;Y)$. Let $(y_i)$ be a sequence of simple functions such that $y_i \rightarrow y$ in $L^1(\Omega;X)$. Then, as in the proof of Lemma~\ref{lem:takeoutwhatisknowni}, $y_i \1_{\{||y_i||_X \leq ||y||_X + 1 \}} \rightarrow y$ almost surely, for all $i = 1, \hdots$. Write $A_i := \{ \omega \in \Omega: ||y_i(\omega)||_X \leq ||y(\omega)||_X + 1\}$. Then $x \cdot y_i\1_{A_i} \rightarrow x \cdot y$ almost surely and $||x \cdot y_i\1_{A_i}||_Y \leq ||x||_X \ ( ||y||_X + 1)$ for all $i = 1, \hdots$. Furthermore $y_i\1_{A_i}$ is a simple function, so
\[ \E [ x \cdot y_i \1_{A_i} |\mathcal G ] = x \E [ y_i \1_{A_i} | \mathcal G].\]

By Lemma~\ref{lem:condexpalmostsureconvergence},
\[ \E \left[ x \cdot y_i \1_{A_i} | \mathcal G \right] \rightarrow \E \left[ x \cdot y | \mathcal G \right]  \quad \textrm{almost surely}, \]
and 
\[ x \cdot \E \left[y_i \1_{A_i} |\mathcal G \right] \rightarrow x \cdot \E \left[ y | \mathcal G \right] \quad \textrm{almost surely}, \]
and the claimed result follows.
\hfill $\square$
\end{quote}

\subsection{Lemma} 
Suppose $x, y \in L^1(\Omega; X)$ such that $x \cdot y \in L^1(\Omega; Y)$ and suppose $x$ and $y$ are independent. Then $\E [ x \cdot y ] = \E[x] \E[y]$.
\label{lem:vectorvaluedindependence}

\begin{quote}
\textsc{Proof:} 
Suppose $x$ is $\mathcal F$-measurable and $y$ is $\mathcal G$-measurable, so $\mathcal F$ and $\mathcal G$ are independent. Then, by applying the tower property (Lemma~\ref{lem:towerproperty}), taking out what is known (Lemma~\ref{lem:takeoutwhatisknownii}) and by using independence (Lemma~\ref{lem:condexpectationindependence}), we find that
\[ \E \left[ x \cdot y \right] = \E \left[ \E \left[ x \cdot y | \mathcal F \right] \right] = \E \left[ x \E \left[ y |\mathcal F \right] \right] = \E \left[ x \E [ y ] \right] = \E [x] \E [y].\]
\hfill $\square$
\end{quote}

\subsection{Lemma (Cauchy-Schwarz)}
\label{lem:cauchyschwartz}
Let $(S, \Sigma, \mu)$ be a $\sigma$-finite measure space.
Let $x, y \in L^2(S;X)$. Then
\[ \left|\left| \int_S x \cdot y  \ d \mu \right|\right|_Y^2 \leq \left|\left| \int_S x^2 \ d \mu \right|\right|_Y \ \left|\left| \int_S y^2 \ d \mu \right|\right|_Y.\]

\begin{quote}
\textsc{Proof:}
Let $\lambda > 0$. Then $x - \lambda y \in L^2(S;X)$, and by Lemma~\ref{lem:integralpositive},
\[ \int_S (x - \lambda y)^2 \ d \mu \succeq 0,\]
so
\[ 2 \lambda \int_S x \cdot y \ d \mu \preceq \int_S x^2 \ d \mu + \lambda^2 \int_S y^2 \ d \mu.\]
Similarly (by considering $(x + \lambda y)^2$,
\[ - 2 \lambda \int_S x \cdot y \ d \mu \preceq \int_S x^2 \ d \mu + \lambda^2 \int_S y^2 \ d \mu,\]
so
\[ 2 \lambda \left| \int_S x \cdot y \ d \mu \right| \preceq \left| \int_S x^2 \ d \mu + \lambda^2 \int_S y^2 \ d \mu \right|, \]
and therefore
\[ \left|\left| \int_S x \cdot y \ d \mu \right|\right|_Y \leq \frac 1 {2 \lambda} \left|\left| \int_S x^2 \ d \mu \right|\right|_Y + \frac \lambda 2 \left|\left| \int_S y^2 \ d \mu \right|\right|_Y,\]
for any $\lambda > 0$.
Now choose \[ \lambda = \sqrt { \frac{\left|\left| \int_S x^2 \ d \mu \right|\right|_Y}{\left|\left| \int_S y^2 \ d \mu \right|\right|_Y} }.\]
\hfill $\square$
\end{quote}


Let $X$ and $Y$ be a Banach spaces, $\cdot : X \times X \rightarrow Y$ a multiplication operation and $(\Omega, \mathcal F, \P)$ a probability space. 
Let $(\mathcal F_t)_{0 \leq t \leq T}$ be a filtration on $\Omega$. As for real-valued stochastic processes, here we say that a process $x : [0,T] \times \Omega \rightarrow X$ is $(\mathcal F_t)$-adapted if for all $t \in [0,T]$ the vector valued random variable $x(t,\cdot)$ is $\mathcal F_t$-measurable.

\subsection{Martingales}
In the remainder of this work, we assume that the stochastic process $(M(t))_{t \geq 0}$ is a cadlag martingale starting at 0.

Under these conditions there exists (see \cite{rogerswilliams}, Theorem VI.36.6) a unique increasing process $[M]$, called the \emph{quadratic variation process of $M$}, such that 
\begin{itemize}
\item[(i)]
$[M](0) = 0$;
\item[(ii)]
$M^2 - [M]$ is a uniformly integrable martingale;
\item[(iii)]
$\Delta [M] = (\Delta M)^2$ on $(0,\infty)$.
\end{itemize}
Here, for a cadlag function $f$,  $\Delta f$ is the function on $(0, \infty)$ with
\[ (\Delta f)(t) := f(t) - f(t-).\]

Since $[M]$ is increasing, it follows that $[M]$ is locally of finite variation, almost surely. Result (iii) tells us that $[M]$ is a cadlag process as well. Because $[M]$ is of finite variation, the Lebesgue-Stieltjes integral 
\[ \int_0^T x(s) \ d [M](s) \]
is defined a.s. for functions $x$ that are integrable with respect to the measure $\mu_{[M]}$ associated to $[M]$.

For an interval $I \subset [0, \infty)$, we write $L^2_M( \Omega \times I; X ) = L^2((\Omega, \P) \times (I,\mu_{[M]}); X)$
for the space of (equivalence classes of) functions $x : \Omega \times I \rightarrow X$, measurable with respect to $\P \otimes \mu_{[M]}$ with
\[ \E \left[ \int_I ||x(s)||_X^2 \ d [M](s) \right] < \infty.\] 

Furthermore we write
$L^{2,a}_M(\Omega \times I;X)$ for the space of adapted functions in $L^2_M(\Omega \times I;X)$.



\subsection{Elementary processes}
We call a process $x \in L^{2,a}_M(\Omega \times [0,T];X)$ \emph{elementary} (denoted $x \in \elem_M$) if it has the form
\begin{equation} \label{eq:elementary} x(t,\omega) = \sum_{i = 0}^n x_i(\omega) \1_{[t_i, t_{i+1})}(t), \quad t \in [0,T], \omega \in \Omega, \end{equation}
where $0 = t_0 \leq t_1 \leq \hdots \leq t_n \leq t_{n+1} = T, $
and where $x_i \in L^2((\Omega,\mathcal F_{t_i}); X)$ for all $i = 0, \hdots, n$.

Because $x \in L^2_M(\Omega \times [0,T];X)$, necessarily
\begin{equation} \label{eq:Mintegrable} \E \left[ ||x_i||_X^2 ([M](t_{i+1}) - [M](t_i)) \right] <\infty, \quad i = 0, \hdots, n.\end{equation}

Note that for $M=W$ a Brownian motion, we have that~(\ref{eq:Mintegrable}) holds for all $x_i \in L^2(\Omega;X)$, $i=1, \hdots, n$, by independence of $\sigma(W(t_{i+1})-W(t_i))$ and $\mathcal F_{t_i}$.

The next result is a slightly modified version of Proposition 3.3 in \cite{vangaans}.
\subsection{Lemma}
\label{lem:shift}
Let $X$ be a Banach space and let $x \in L^2_M(\Omega \times [0,T]; X)$. Define, for $t \in (0,T)$,
\[ x_t(s) = \left\{ \begin{array}{cl} x(s-t), \quad & s \in [t,T] \\
0, \quad & s \in [0,t). \end{array} \right.
\]
Then
\[ \lim_{t \downarrow 0} ||x_t - x||_{L^2_M(\Omega \times [0,T];X)} = 0.\]
\begin{quote}
\textsc{Proof:}
If $x$ has the form $x = a \1_{A \times I}$, where  $A \in \mathcal F$, $I \in \mathcal B([0,T])$ and $a \in X$, then
\[ ||x_t(s) - x(s)||_X \rightarrow 0 \quad \textrm{as} \ t \downarrow 0 \quad \textrm{for almost all} \ (\omega,s) \in [0,T] \times \Omega.\]
By dominated convergence,
\[ \E \left[ \int_0^T ||x_t(s) - x(s)||_X^2 \ d [M](s) \right] \rightarrow 0 \quad \textrm{as} \ t \downarrow 0.\]
Therefore the result holds when $x$ is a simple function.
The general case follows by approximation.
\hfill $\square$
\end{quote}

\subsection{Lemma}
\label{lem:elemdenseinL2}
$\elem_M$ is dense in $L^{2,a}_M$.

\begin{quote}
\textsc{Proof:}
Let $x \in L^{2,a}_M(\Omega \times [0,T]; X)$.
Let $\epsilon > 0$ be given. By Lemma~\ref{lem:shift}, we can choose $N \geq 2$ so large that
\[ || x - x_{\frac T N} ||_{L^2_M(\Omega \times [0,T];X)} < \half \epsilon.\]
Because $x$ is adapted, we see that $x_{\frac T N}(t)$ is $\mathcal F_{t-\frac T N}$-measurable for all $t \in [0,T]$.
For $k = 1, \hdots, N-1$, since $x_{\frac T N}(t)$ is $\mathcal F_{\frac{kT}{N}}$-measurable for all $t \in \left[ \frac{kT}{N}, \frac{(k+1)T}{N} \right)$, the restricted function $\left. x_{\frac T N} \right|_{\left[\frac {kT}{N}, \frac {(k+1)T}{N} \right) }$ can be approximated by a simple function $\phi_k : \left[ \frac{kT}{N}, \frac{(k+1)T}{N} \right) \rightarrow L^2((\Omega; \mathcal F_{\frac{kT}{N}}), X)$, such that $\phi_k \in L^2_M(\Omega \times[\frac {kT} N, \frac{(k+1)T} N ); X)$ and
\[ \left|\left| \left. x_{\frac T N} \right|_{\left[\frac {kT}{N}, \frac {(k+1)T}{N} \right) } - \phi_k \right|\right|_{L^2_M \left(\Omega \times \left[ \frac {kT}{N}, \frac{(k+1)T}{N} \right); X \right)} < \frac \epsilon { 2 (N-1)}.\]

Define $\phi: [0,T] \rightarrow L^2(\Omega;X)$ by gluing together the functions $\phi_k$:
\[ \phi(t) = \left\{ \begin{array}{cl} \phi_k(t), \quad & t \in \left[ \frac {kT}{N}, \frac{(k+1)T}{N} \right), \quad k = 1, \hdots, N-1, \\ \\
0, \quad & t \in \left[ 0, \frac T N \right) \end{array} \right. \]

Then, for $t \in \left[ \frac {kT}{N} , \frac {(k+1)T}{N} \right)$, because $\phi(t)$ is $\mathcal F_{\frac {kT}{N}}$-measurable, $\phi(t)$ is $\mathcal F_t$-measurable, so $\phi$ is adapted, i.e. $\phi \in \elem_M$.

Furthermore,
\[ || x_{\frac T N} - \phi ||_{L^2_M(\Omega \times [0,T]; X)} < (N-1) \frac \epsilon {2 (N-1)} = \half \epsilon,\]
which shows that
\[ || x - \phi ||_{L^2_M(\Omega \times [0,T];X)} < \epsilon.\]
\hfill $\square$
\end{quote}

\subsection{Definition (Stochastic integral of an elementary process)}
For $x \in \elem_M$, $x = \sum_{i=0}^n x_i \1_{[t_i, t_{i+1})}$ define the \emph{integral of $x$ with respect to $M$} by
\[ \int_0^t x(s) \ d M(s) = \sum_{i= 0}^n x_i \left(M(t_{i+1} \wedge t) - M(t_i \wedge t) \right).\]




\subsection{Lemma (The It\^o isometry)}
\label{lem:itoisom}
If $x \in \elem_M$ then
\begin{equation} \label{eq:itoisom} \E \left[ \left(\int_0^t x(s) \ d M(s) \right)^2 \right] = \E \left[ \int_0^t x(s)^2 \ d [M](s) \right], \quad t \in [0,T].\end{equation}
As a consequence, $\int_0^t x(s) \ d M(s) \in L^2_M(\Omega \times [0,T];X)$.

\begin{quote}
\textsc{Proof:} It suffices to show this for $t = T$, because for $t < T$ the function $s \mapsto x(s) \1_{s \leq t}$ is elementary if $x$ is elementary.

Write $x = \sum_{i=0}^n x_i \1_{[t_i, t_{i+1})}$. For $i,j = 0, \hdots, n$, if $i < j$, by Lemma~\ref{lem:towerproperty} and Lemma~\ref{lem:takeoutwhatisknownii},
\begin{align*} & \E \left[(x_i \cdot x_j) (M(t_{i+1}) - M(t_i))  (M(t_{j+1}) - M(t_j) ) \right] \\
& = \E \left[ \E [ (x_i \cdot x_j) (M(t_{i+1}) - M(t_i)) (M(t_{j+1}) - M(t_j)) | \mathcal F_{t_j} ] \right] \\
& = \E \left[ (x_i \cdot x_j) (M(t_{i+1}) - M(t_i)) \E [ M(t_{j+1}) - M(t_j) | \mathcal F_{t_j} ] \right] = 0. \\
\end{align*}

If $i = j$, then by Lemma~\ref{lem:takeoutwhatisknownii}, and by the definition of the quadratic variation $[M]$,
\begin{align*}
\E \left[x_i^2 (M(t_{i+1}) - M(t_i))^2 \right] & = \E \left[ \E [ x_i^2 (M(t_{i+1})-M(t_i))^2 |\mathcal F_{t_i} ] \right] \\
& = \E \left[ x_i^2 \E [ (M(t_{i+1}) - M(t_i))^2 | \mathcal F_{t_i} ] \right] \\
& = \E \left[ x_i^2 \E [ M(t_{i+1})^2 - M(t_i)^2 | \mathcal F_{t_i} ] \right] \\
& = \E \left[ x_i^2 \E [ [M](t_{i+1}) - [M](t_i) | \mathcal F_{t_i} ] \right] \\
& = \E \left[ x_i^2 ( [M](t_{i+1}) - [M](t_i) ) \right].
\end{align*}
Hence
\begin{align*}
\E \left[ \left( \int_0^T x(s) \ d M(s) \right)^2 \right] & = \sum_{i,j = 0}^n \E \left[ (x_i \cdot x_j) (M(t_{i+1}) - M(t_i)) (M(t_{j+1})-M(t_j)) \right] \\
& = \sum_{i=0}^n \E \left[ x_i^2 ([M](t_{i+1}) - [M](t_i)) \right] \\
& = \E \left[ \int_0^T x^2(s) \ d [M](s) \right].
\end{align*}
\hfill $\square$
\end{quote}

As in the real-valued case, the It\^o isometry will enable us to define the stochastic integral for a broader class of integrands. The first step is to define the class of integrable functions. Loosely speaking, we will construct a space $\mathcal M_M$ of functions $x: \Omega \times [0,T] \rightarrow X$ that are integrable with respect to $M$. Furthermore, we construct a space $\mathcal M_t$ so that $\int_0^t x(s) \ d M(s) \in \mathcal M_t$, so the evaluation of the stochastic integral at time $t$ belongs to $\mathcal M_t$, for $x \in \mathcal M_M$ and $t \in [0,T]$.

\subsection{Lemma} 
\label{lem:mnorm}
For any measure space $(S, \mu)$ the mapping
\[ ||\cdot||_{\mathcal M(S;X)} : L^2(S;X) \rightarrow \reals\]
defined by
\[ ||x||_{\mathcal M(S;X)} =  \left|\left| \int_S x^2 \ d \mu \right|\right|_Y^{\half}\]
is continuous and defines a norm on $L^2(S;X)$. Furthermore the important estimate
\begin{equation} \label{eq:estimateMnorm} ||x||_{\mathcal M(S;X)} \leq ||x||_{L^2(S;X)} \end{equation}
holds.
\label{lem:Mnorm}

\begin{quote}
\textsc{Proof:}
Because
\[ || x||_{\mathcal M(S;X)}^2 = \left|\left| \int_S x^2 \ d \mu \right|\right|_Y \leq \int_S ||x||^2_X \ d \mu,\]
if $x \in L^2(S;X)$, then $||x||_{\mathcal M(S;X)} < \infty$, and if $x = 0$, then $||x||_{\mathcal M(S;X)} = 0$.

Obviously, for $\alpha \in \reals$ and $x \in L^2(S;X)$
\[ ||\alpha x||_{\mathcal M(S;X)} = |\alpha| \ ||x||_{\mathcal M(S;X)}.\]

Let $x, y \in L^2(S;X)$. Then
\begin{align*}
\left|\left| x + y \right|\right|^2_{\mathcal M(S;X)} & = \left|\left| \int_S (x + y)^2 \ d \mu \right|\right|_Y \\
& \leq \left|\left| \int_S x^2 \ d \mu \right|\right|_Y + \left|\left| \int_S y^2 \ d \mu \right|\right|_Y + 2 \left|\left| \int_S x \cdot y \ d \mu \right|\right|_Y.
\end{align*}
By Lemma~\ref{lem:cauchyschwartz},
\[ \left|\left| \int_S x \cdot y \ d \mu \right|\right|_Y \leq \left|\left| \int_S x^2 \ d \mu \right|\right|_Y^{\half} \ \left|\left| \int_S y^2 \ d \mu \right|\right|_Y^{\half},\]
which shows that
\[ \left|\left| x + y \right|\right|^2_{\mathcal M(S;X)} \leq \left( \left|\left| \int_S x^2 \ d \mu \right|\right|_Y^{\half} + \left|\left| \int_S y^2 \ d \mu \right|\right|_Y^{\half} \right)^2 = \left( ||x|| + ||y|| \right)^2.\]

Suppose now that $x \neq 0$ in $L^2(S;X)$. Therefore $x \neq 0$ in $X$ on a set $A$ with $\mu(A) > 0$. Then, by assumption (iii) in Definition~\ref{def:multiplication}, $x^2 \neq 0$ on $A$. Lemma~\ref{lem:integralzero} shows that 
\[ \int_S x^2 \1_A \ d \mu \neq 0,\]
and therefore
\[ \left|\left| \int_S x^2 \ d \mu \right|\right|_Y \geq \left|\left| \int_S x^2 \1_A \ d \mu \right|\right|_Y > 0.\]
\hfill $\square$
\end{quote}

\subsection{Definition (Stochastically integrable functions)}
\label{def:stochintegrable}
With $||\cdot||_{\mathcal M(S;X)}$ as in Lemma~\ref{lem:Mnorm}, define $\mathcal M(S;X)$ to be the completion of $L^2(S;X)$ with respect to $||\cdot||_{\mathcal M(S;X)}$.

Let $||\cdot||_{\mathcal M_M} := ||\cdot||_{\mathcal M_M(\Omega \times [0,T];X)} := ||\cdot||_{\mathcal M(\Omega\times [0,T];X)}$ as defined in Lemma~\ref{lem:Mnorm} where we take $\P \times \mu_{[M]}$ as the measure on $\Omega \times [0,T]$.
Define $\mathcal M_M$, the \emph{space of stochastically integrable functions with respect to $M$}, to be the completion of $\elem_M$ with respect to $||\cdot||_{\mathcal M_M}$.

Furthermore define $\mathcal M_t$ to be the completion of $L^2((\Omega, \mathcal F_t);X)$ with respect to $||\cdot||_{\mathcal M(\Omega;X)}$.

\subsection{Lemma}
\label{lem:L2denseinR} We have
\[ \mathcal M_M = \overline { \elem_M } \subset \overline { L^{2,a}_M(\Omega \times [0,T]; X)} = \mathcal M_M,\]
where the completion is taken with respect to the $\mathcal M_M$-norm.

\begin{quote}
\textsc{Proof:} Let $x \in L^2_M(\Omega\times[0,T];X))$ be adapted. By Lemma~\ref{lem:elemdenseinL2} there exists a sequence of elementary functions $x_n \in \elem_M$ such that
\[ ||x_n - x||_{L^2_M(\Omega \times [0,T];X)} \rightarrow 0, \quad \textrm{as} \ n \rightarrow \infty.\]
But
\begin{align*} ||x_n - x||^2_{\mathcal M_M} & = \left|\left| \E \left[ \int_0^T (x_n-x)^2(s) \ d[M](s) \right] \right|\right|_Y \\
& \leq \E \left[ \int_0^T \left|\left| (x_n - x)^2 (s) \right|\right|_Y \ d[M](s) \right] \\
& \leq \E \left[ \int_0^T \left|\left| x_n(s) - x(s) \right|\right|^2_X \ d[M](s) \right] = ||x_n - x||_{L^2_M(\Omega \times [0,T];X)}^2,
\end{align*}
which shows that $x \in \mathcal M_M$. Because $\elem_M \subset L^{2,a}_M(\Omega \times [0,T]; X)$, the claimed result follows.
\hfill $\square$
\end{quote}

\subsection{Lemma}
If $X$ is a Hilbert space, and we use the standard construction of $Y$ as described in Example~\ref{multiplicationexamples}~(iv),
then 
\[ ||\cdot||_{\mathcal M_M} = ||\cdot||_{L^2_M(\Omega \times [0,T]; X)} \quad \textrm{and} \quad \mathcal M_M = L^{2,a}_M(\Omega \times [0,T];X).\]
\begin{quote}
\textsc{Proof:}
In this case, for $x \in \mathcal M_M$,
\begin{align*} ||x||_{\mathcal M_M}^2 & = \left|\left| \E \left[ \int_0^t x^2(s) \ d[M](s) \right] \right|\right|_Y = \left| \E \left[ \int_0^t (x(s),x(s))_X \ d[M](s) \right] \right| \\
& = \E \left[ \int_0^t ||x(s)||_X^2 \ d[M](s) \right] = ||x||_{L^2_M(\Omega \times [0,T];X)}.
\end{align*}
\hfill $\square$
\end{quote}

\subsection{Theorem and Definition (It\^o integral)}
\label{thm:itointegral}
Let $x \in \mathcal M_M$. For $t \in [0,T]$ there exists a unique $z_t \in \mathcal M_t$ such that for any sequence of elementary functions $(x_n)_{n \in \ints}$ in $\elem_M$ with $||x_n - x||_{\mathcal M_M} \rightarrow 0$ as $n \rightarrow \infty$,
\[ \int_0^t x_n(s) \ d M (s) \rightarrow z_t \quad \textrm{in} \ \mathcal M_t.\]
This $z_t$ will be called the \emph{It\^o integral} of $x$ with respect to $M$ and will be denoted by
\[ \int_0^t x(s) \ d M(s).\]

\begin{quote}
\textsc{Proof:}
Let $x \in \mathcal M_M$ and suppose $(x_n)_{n \in \ints}$ is a sequence in $\elem_M$ such that $x_n \rightarrow x$ in $\mathcal M_M$. 
For $t \in [0,T]$, by the It\^o isometry (Lemma~\ref{lem:itoisom}), $\int_0^t x_n(s) \ d M(s) \in L^2((\Omega, \mathcal F_t); X)$, and
\begin{align*}
\left|\left| \int_0^t x_n(s) \ d M(s) - \int_0^t x_m(s) \ d M(s) \right|\right|_{\mathcal M_t}^2 & =  \left|\left| \E\left[ \left(\int_0^t (x_n - x_m)(s) \ d M(s) \right)^2 \right] \right|\right|_Y \\
& = \left|\left| \E \left[ \int_0^t (x_n - x_m)^2(s) \ d[M](s) \right] \right|\right|_Y.
\end{align*}
Because $(x_n - x_m)^2(s) \succeq 0$ in $Y$ for almost all $s \in [0,T]$, almost surely, by Lemma~\ref{lem:integralpositive},
\[ \E \left[ \int_0^t (x_n - x_m)^2(s) \ d[M](s) \right] \preceq \E \left[ \int_0^T (x_n - x_m)^2(s) \ d[M](s) \right],\]
and hence
\begin{align*} \left|\left| \int_0^t x_n(s) - x_m(s) \ d M(s) \right|\right|_{M_t} & = \left|\left| \E \left[ \int_0^t (x_n - x_m)^2(s) \ d[M](s) \right] \right|\right|_Y \\
& \leq \left|\left| \E \left[ \int_0^T (x_n - x_m)^2(s) \ d[M](s) \right] \right|\right|_Y \\
& = ||x_n - x_m||_{\mathcal M_M}^2.
\end{align*}

This shows that $\left(\int_0^t x_n \ d M(s)\right)_{n \in \ints}$ is a Cauchy sequence in $\mathcal M_t$. Since $\mathcal M_t$ is, by definition, complete, this sequence has a limit $z_t \in \mathcal M_t$.

Suppose that, besides $z_t$, also $\zeta_t \in \mathcal M_t$ such that $\int_0^t \xi_n(s) \ d M(s) \rightarrow \zeta_t$ for some sequence $\xi_n \rightarrow x$ in $\mathcal M$. Then
\begin{align*}
|| z_t - \zeta_t||_{\mathcal M_t}^2 & = \lim_{n \rightarrow \infty} \left|\left| \int_0^t x_n(s) \ d M(s) - \int_0^t \xi_n(s) \ d M(s) \right|\right|_{\mathcal M_t}^2 \\
& = \lim_{n \rightarrow \infty} \left|\left| \E \left[ \left( \int_0^t x_n(s) - \xi_n(s) \ d M(s) \right)^2 \right] \right|\right|_Y \\
& = \lim_{n \rightarrow \infty} \left|\left| \E \left[ \int_0^t (x_n(s) - \xi_n(s))^2 \ d[M](s) \right] \right|\right|_Y \\
& \leq \lim_{n \rightarrow \infty} \left|\left| \E \left[ \int_0^T (x_n(s) - \xi_n(s))^2 \ d[M](s) \right] \right|\right|_Y \\
& = \lim_{n \rightarrow \infty} \left|\left|x_n - \xi_n \right|\right|_{\mathcal M_M}^2 \\
& \leq \lim_{n \rightarrow \infty} (||x_n - x ||_{\mathcal M_M} + ||x - \xi_n||_{\mathcal M_M})^2 = 0,
\end{align*} 
which shows that $z_t = \zeta_t$ and hence that the integral is uniquely defined.
\hfill $\square$
\end{quote}

\section{Properties of $\mathcal M(S;X)$}
In Section~\ref{sec:stochasticintegral} we introduced the space $\mathcal M(S;X)$. In this section some properties of this space will be developed. Throughout this section the spaces $X, Y$ and the mapping $(x_1,x_2) \mapsto x_1 \cdot x_2 : X \times X \rightarrow Y$ satisfy the properties of Definition~\ref{def:multiplication}, and we use the $\mathcal M(\cdot,\cdot)$-spaces as introduced in Definition~\ref{def:stochintegrable}. Furthermore, $(S,\Sigma, \mu)$ is a finite measure spaces and $(\Omega, \mathcal F, \P)$ a probability space.

\subsection{Lemma}
\label{lem:mappings}
The following estimates hold:
\begin{itemize}
\item[(i)]
\[ \left|\left| \E \left[ x | \mathcal G \right] \right|\right|_{\mathcal M(\Omega;X)} \leq ||x||_{\mathcal M(\Omega;X)}.\]
where $x \in L^2((\Omega,\mathcal F);X)$, and $\mathcal G$ a sub-$\sigma$-algebra of $\mathcal F$,
\item[(ii)] 
\[ \left|\left| \int_S A x \ d \mu \right|\right|_{\overline V} = \left|\left| \left( \int_S A x \ d \mu \right)^2 \right|\right|_W^{\half} \leq \sqrt{\mu(S) \ ||B||_{\mathcal L(Y;W)}} \ ||x||_{\mathcal M(S;X)},\]
where
\begin{itemize}
\item $x \in L^2(S;X)$,
\item $V$ is a Banach space, $W$ a Banach lattice, and $(v_1, v_2) \mapsto v_1 \cdot v_2 : V \times V \rightarrow W$ a multiplication that satisfies the conditions of Definition~\ref{def:multiplication} (e.g. take $V=X$, $W=Y$);
\item
$A \in \mathcal L(X;V)$ such that there exists a $B \in \mathcal L(Y;W)$ such that $(A x)^2 \preceq B(x^2)$ for all $x \in X$;
\item
$\overline V$ the completion of $V$ with respect to the norm $||v||_{\overline V} = ||v^2||_W^{\half}$.
\end{itemize}
\item[(iii)]
\[ \left|\left| A x \right|\right|_{\overline V} = \left|\left| (A x)^2 \right|\right|_W^{\half} \leq \sqrt{ ||B||_{\mathcal L(Y;W)} } \ ||x||_X,\]
with the same notation and conditions as under (ii).
\item[(iv)]
\[ \left|\left| \int_S A x^2 \ d \mu \right|\right|_Z \leq ||A||_{\mathcal L(Y;Z)} \ ||x||_{\mathcal M(S;X)}^2,\]
where $A \in \mathcal L(Y;Z)$ with $Z$ a normed space, and $x \in L^2(S;X)$.
\item[(v)]
\[ \left|\left| \int_S x \cdot y \ d \mu \right|\right| \leq ||x||_{\mathcal M(S;X)} \ ||y||_{\mathcal M(S;X)},\]
where $x,y \in L^2(S;X)$.
\end{itemize}

\begin{quote}
\textsc{Proof:}
\item[(i)] By Jensen's inequality for vector valued conditional expectation (Theorem~\ref{thm:cjensen}),
\begin{align*} \left|\left| \E \left[ x | \mathcal G \right] \right|\right|_{\mathcal M(\Omega;X)}^2 & = \left|\left| \E \left[ \left(\E \left [ x  | \mathcal G \right]\right)^2 \right] \right|\right|_Y \leq \left|\left| \E \left[ \E \left[ x^2 |\mathcal G \right] \right] \right|\right|_Y \\
& = \left|\left| \E \left[ x^2 \right] \right|\right|_Y = ||x||_{\mathcal M(\Omega;X)}^2.
\end{align*}
\item[(ii)]
Again by using Jensen's theorem (this time Corollary~\ref{cor:jensen}), and by Lemma~\ref{lem:integralpositive}, the result follows from the estimate
\begin{align*}
\left|\left| \int_S A x \ d \mu \right|\right|_{\overline V}^2 & = \left|\left| \left(\int_S A x \ d \mu \right)^2 \right|\right|_W \leq \left|\left| \mu(S) \int_S (A x)^2 \ d \mu \right|\right|_W \\
& \leq |\mu(S)| \left|\left| \int_S B x^2 \ d \mu \right|\right|_W = |\mu(S)| \ \left|\left| B \int_S x^2 \ d \mu \right|\right|_W \\
& \leq |\mu(S)| \ ||B ||_{\mathcal L(Y;W)} \ \left|\left| \int_S x^2 \ d \mu \right|\right|_Y = |\mu(S)| \ ||B||_{\mathcal L(Y;W)} \ ||x||_{\mathcal M(S;X)}^2.
\end{align*}
\item[(iii)]
This is an immediate consequence of (ii) by letting $S$ consist of only one point, having measure one.
\item[(iv)]
\begin{align*} \left|\left| \int_S A x^2 \ d \mu \right|\right|_Z & = \left|\left| A \int_S x^2 \ d \mu \right|\right|_Z \leq ||A||_{\mathcal L(Y;Z)} \ \left|\left| \int_S x^2 \ d \mu \right|\right|_Y \\
& = ||A||_{\mathcal L(Y;Z)} \ ||x||_{\mathcal M(S;X)}^2.
\end{align*}
\item[(v)]
A restatement of Cauchy-Schwarz, Lemma~\ref{lem:cauchyschwartz}.
\hfill $\square$
\end{quote}

\subsection{Corollary}
The mappings mentioned in Lemma~\ref{lem:mappings}, i.e. (notation and conditions as in the previous lemma)
\begin{itemize}
\item[(i)]
$\E[\cdot|\mathcal G] : L^2((\Omega,\mathcal F);X) \rightarrow L^2((\Omega, \mathcal G);X)$;
\item[(ii,a)]
$x \mapsto \int_S A x \ d \mu: L^2(S;X) \rightarrow \overline V$;
\item[(ii,b)]
$x \mapsto \left(\int_S A x \ d \mu \right)^2 : L^2(S;X) \rightarrow W$;
\item[(iii,a)]
$x \mapsto A x: X \rightarrow \overline V$;
\item[(iii,b)]
$x \mapsto (A x)^2: X \rightarrow W$;
\item[(iv)]
$x \mapsto \int_S A x^2 \ d \mu: L^2(S;X) \rightarrow Z$;
\item[(v)]
$(x,y) \mapsto \int_S x \cdot y \ d \mu: L^2(S;X) \times L^2(S;X) \rightarrow Y$;
\end{itemize} 
are continuous with respect to the $\mathcal M(S;X)$ norm\footnote{Under (v), with respect to the product topology on $L^2(S;X) \times L^2(S;X)$ induced by the $\mathcal M(S;X)$ norm}. Therefore these mappings are welldefined on $\mathcal M(S;X)$ and the estimates of the lemma remain valid.

\subsection{Lemma}
\label{lem:cartesianproduct}
Let $(X^1, Y^1)$ and $(X^2,Y^2)$ be two pairs of Banach spaces for which the conditions of Definition~\ref{def:multiplication} hold.

Then $X^1 \times X^2$ satisfies the conditions of Definition~\ref{def:multiplication} as well, and
\[ \mathcal M(S;X^1) \times \mathcal M(S;X^2) \cong \mathcal M(S;X^1 \times X^2).\]

\begin{quote}
\textsc{Proof:}
Define multiplication by $\cdot : (X^1 \times X^2) \times (X^1 \times X^2) \rightarrow Y^1 \times Y^2$ by 
\[ (x_1, x_2) \cdot (y_1, y_2) = (x_1 \times y_1, x_2 \times y_2).\]
Then this multiplication satisfies the conditions of Definition~\ref{def:multiplication}.

We find that
\begin{align*}
\left|\left| (x_1, x_2) \right|\right|_{\mathcal M(S;X^1 \times X^2)}^2 & = \left|\left| \int_S (x_1^2(t), x_2^2(t)) \ d \mu \right|\right|_{Y^1 \times Y^2} 
\\
& = \left|\left| \int_S x_1^2 \ d\mu \right|\right|_{Y^1} + \left|\left| \int_S x_1^2 \ d\mu \right|\right|_{Y^2} \\
& = ||x_1||_{\mathcal M(S;X^1)}^2 + ||x_2||_{\mathcal M(S;X^2)}^2 \\
& \leq \left( ||x_1||_{\mathcal M(S;X^1)} + ||x_2||_{\mathcal M(S;X^2)} \right)^2 \\
& = ||(x_1, x_2)||_{\mathcal M(S;X^1) \times \mathcal M(S;X^2)}^2,
\end{align*}
and
\begin{align*} ||(x_1, x_2)||_{\mathcal M(S;X^1) \times \mathcal M(S;X^2)}^2 & = \left( ||x_1||_{\mathcal M(S;X^1)} + ||x_2||_{\mathcal M(S;X^2)} \right)^2 \\
& \leq 2 \left( ||x_1||_{\mathcal M(S;X^1)}^2 + ||x_2||_{\mathcal M(S;X^2)}^2 \right) \\
& = 2 || (x_1, x_2) ||_{\mathcal M(S;X^1 \times X^2)}^2,
\end{align*}
showing that the norms of $\mathcal M(S;X^1 \times X^2)$ and $\mathcal M(S;X^1) \times \mathcal M(S;X^1)$ are equivalent.

Finally, 
\begin{align*} \mathcal M(S;X^1 \times X^2) & = \overline{ L^2(S; X^1 \times X^2)) } = \overline{ L^2(S; X^1) \times L^2(S; X^2)} \\
& = \overline{ L^2(S; X^1) } \times \overline{ L^2(S; X^2) } = \mathcal M(S;X^1) \times \mathcal M(S;X^2),
\end{align*}
which shows that $\mathcal M(S;X^1 \times X^2)$ and $\mathcal M(S;X^1) \times \mathcal M(S;X^2)$ contain the same elements.
\hfill $\square$
\end{quote}

\subsection{Lemma}
Define a bilinear, symmetric mapping $\circ: \mathcal M(S;X) \times \mathcal M(S;X) \rightarrow Y$ by
\[ x_1 \circ x_2 = \int_S x_1 \cdot x_2 \ d \mu.\]
Then this mapping satisfies the requirements of a multiplication as in Definition~\ref{def:multiplication}.
\begin{quote}
\textsc{Proof:}
Note that by Cauchy-Schwarz (Lemma~\ref{lem:cauchyschwartz}),
\[ \left|\left|\int_S x_1 \cdot x_2 \ d \mu \right|\right|_Y \leq \left|\left| \int_S x_1^2 \ d \mu \right|\right|_Y^{\half} \ \left|\left| \int_S x_1^2 \ d \mu \right|\right|_Y^{\half},\]
which makes the mapping `$\circ$' welldefined, and ensures that the multiplication satisfies the first condition, namely
\[ || x \circ y ||_Y \leq ||x||_{\mathcal M(S;X)} \ ||y||_{\mathcal M(S;X)}.\]

Lemmas~\ref{lem:integralpositive} and~\ref{lem:integralzero} show that the other two conditions are also satisfied.
\hfill $\square$
\end{quote}

The above lemma makes it possible for us to discuss the space $\mathcal M(S;\mathcal M(T;X))$, where $(S,\mu)$ and $(T,\lambda)$ are finite measure spaces.

\subsection{Lemma}
For $S$ and $T$ finite measure spaces, and $(X,Y)$ as usual, there exists an isometric isomorphism betwenen $\mathcal M(S;\mathcal M(T;X))$ and $\mathcal M(S \times T; X)$.

\begin{quote}
\textsc{Proof:}
First we will show that the norms are equivalent.
\footnote{This is the essence of the proof; the rest are verifications to be mathematically precise.}
\begin{align}
\label{eq:mnormsequivalent}
||x||_{\mathcal M(S;\mathcal M(T;X))}^2 & = \left|\left| \int_S x \circ x \ d \mu \right|\right|_Y \\
\nonumber & = \left|\left| \int_S \left(\int_T x^2 \ d \lambda\right) \ d \mu \right|\right|_Y \\
\nonumber & = \left|\left| \int_{S \times T} x^2 \ \left( d \mu \otimes d \lambda\right) \right|\right|_Y \\
\nonumber & = ||x||_{\mathcal M(S \times T;X)}^2.
\end{align}

Next we show that the two spaces contain `the same elements'. 
We have
\[L^2(S \times T;X) \cong L^2(S; L^2(T;X)) \hookrightarrow L^2(S; \mathcal M(T;X)) \hookrightarrow \mathcal M(S; \mathcal M(T;X)),\]
so that
\[ \mathcal M(S \times T;X) = \overline{L^2(S \times T;X)} \hookrightarrow \overline{ \mathcal M(S;\mathcal M(T;X))} = \mathcal M(S;\mathcal M(T;X))\]
since the norms of $\mathcal M(S \times T;X)$ and $\mathcal M(S;\mathcal M(T;X))$ are equal, as shown above.

Write $\iota: L^2(T;X) \hookrightarrow \mathcal M(T;X)$, $\gamma: L^2(S; \mathcal M(T;X)) \hookrightarrow \mathcal M(S; \mathcal M(T;X))$ and $\eta: L^2(S \times T;X) \hookrightarrow L^2(S;\mathcal M(T;X))$ for some linear embeddings that we will need. Note that, by Lemma~\ref{lem:mnorm}, $||\gamma||\leq 1$, $||\iota||\leq 1$ and $||\eta|| \leq 1$.

Let $\epsilon > 0$. Now if $x \in \mathcal M(S;\mathcal M(T;X))$, then there exists a $y \in L^2(S; \mathcal M(T;X))$ such that $||x - \gamma(y) ||_{\mathcal M(S;\mathcal M(T;X))} < \epsilon / 3$. In turn, $y$ can be approximated by a simple function $z = \sum_{i=1}^n z_i \1_{A_i}$ where the $A_i$ are disjoint, measurable sets, and $z_i \in \mathcal M(T;X)$ for all $i = 1, \hdots, n$, such that $||y - z||_{L^2(S;\mathcal M(T;X))} < \epsilon / 3$. For $i = 1, \hdots, n$, choose $\zeta_i \in L^2(T;X)$ so that $||\iota(\zeta_i) - z_i ||_{\mathcal M(T;X)} < \epsilon / (3 n \sqrt{\mu(A_i)})$. Then $\zeta := \sum_{i=1}^n \zeta_i \1_{A_i} \in L^2(S;L^2(T;X)) \cong L^2(S \times T;X)$, and 
\begin{align*} 
||\eta(\zeta) - z ||_{L^2(S;\mathcal M(T;X))} & = \left|\left| \sum_{i=1}^n (\iota(\zeta_i) - z_i) \1_{A_i} \right|\right|_{L^2(S;\mathcal M(T;X))} \\
& \leq \sum_{i=1}^n \left|\left| (\iota(\zeta_i) - z_i) \1_{A_i} \right|\right|_{L^2(S;\mathcal M(T;X))} \\
& = \sum_{i=1}^n \sqrt{\mu(A_i)} ||\iota(\zeta_i) - z_i||_{\mathcal M(T;X)} < \frac \epsilon 3.
\end{align*}

We conclude that
\begin{align*} || (\gamma \circ \eta)(\zeta) - x ||_{\mathcal M(S;\mathcal M(T;X))} & \leq ||(\gamma \circ \eta)(\zeta) - \gamma(z) ||_{\mathcal M(S;\mathcal M(T;X)} + ||\gamma(z - y)||_{\mathcal M(S;\mathcal M(T;X))} \\
& \quad + ||\gamma(y) - x||_{\mathcal M(S;\mathcal M(T;X))} \\
& \leq ||(\gamma \circ \eta)(\zeta) - z||_{L^2(S;\mathcal M(T;X))} + ||z - y||_{L^2(S;\mathcal M(T;X))} \\
& \quad + ||\gamma(y) - x||_{\mathcal M(S;\mathcal M(T;X))} \\
& < \epsilon,
\end{align*}
which shows that $L^2(S \times T;X)$ is densely embedded in $\mathcal M(S;\mathcal M(T;X))$, i.e. $\overline{(\gamma \circ \eta) (L^2(S \times T;X))} = \mathcal M(S; \mathcal M(T;X))$. 

Now because for $x \in L^2(S \times T;X)$, we have that $||(\gamma \circ \eta)(x)||_{\mathcal M(S;\mathcal M(T;X))} = ||x||_{\mathcal M(S\times T;X)}$ as shown by~(\ref{eq:mnormsequivalent}), we find that
\begin{align*}
\mathcal M(S; \mathcal M(T;X)) & = \overline{(\gamma \circ \eta) (L^2(S \times T;X))} = (\gamma \circ \eta) \left(\overline {L^2(S \times T;X)} \right) \\
& = (\gamma \circ \eta) \left( \mathcal M(S \times T;X) \right), 
\end{align*}
which gives us the isometric isomorphism between $\mathcal M(S \times T;X)$ and $\mathcal M(S; \mathcal M(T;X))$.
\hfill $\square$
\end{quote}

\section{Properties of the stochastic integral}
In this section we will verify several properties of stochastic integration, as developed in Section~\ref{sec:stochasticintegral}.
Unless stated otherwise, $X$ is a Banach space with multiplication $\cdot: X \times X \rightarrow Y$ satisfying the conditions of Definition~\ref{def:multiplication}, and $M$ is a cadlag martingale starting at 0.

\subsection{Proposition (It\^o's isometry)}
Let $x \in \mathcal M_M(\Omega \times [0,T];X)$. Then for $t \in [0,T]$, both
\[ \E \left[ \left( \int_0^t x(s) \ d M(s) \right)^2 \right] \quad \textrm{and} \quad \E \left[ \int_0^t x^2(s) \ d[M](s) \right] \] are welldefined, and equal.

\begin{quote}
\textsc{Proof:}
Let $(x_n)_{n\in \ints}$ be an approximation of $x$ by elementary functions in $\elem_M$. 
Then, by Theorem~\ref{thm:itointegral},
\[ \lim_{n \rightarrow \infty} \E \left[\left(\int_0^t x_n(s) \ d M(s) \right)^2 \right] =:\E \left[ \left(\int_0^t x(s) \ d M(s) \right)^2 \right],\]
and
\[ \lim_{n \rightarrow \infty} \int_0^t \E \left[ x_n^2 \right] \ d[M](s) =: \int_0^t \E \left[ x^2(s) \right] \ d[M](s)\]
are welldefined elements of $Y$.

Because, for all $n \in \ints$,
\[ \E \left[ \left( \int_0^t x_n(s) \ d M(s) \right)^2 \right] = \E \left[ \int_0^t x_n(s)^2 \ d[M](s) \right],\]
the limits are equal.
\hfill $\square$
\end{quote}

\subsection{Lemma}
The mapping $t \mapsto \int_0^t x(s) \ d M(s)$ is an element of 
\[ L^{2,a}([0,T];\mathcal M_T) := \overline{ \{ x \in L^2([0,T];L^2(\Omega;X)) : x \ \textrm{adapted} \} },\] 
where the closure is taken with respect to the $L^2([0,T];\mathcal M_T)$-norm.

\begin{quote}
By the It\^o isometry, for $x \in \mathcal M_M$,
\begin{align*}
\int_0^T \left|\left|\int_0^t x(s) \ d M(s) \right|\right|_{\mathcal M_T}^2 \ d t & = \int_0^T \left|\left| \E \left[ \left( \int_0^t x(s) \ d M(s) \right)^2 \right] \right|\right|_Y \ d t \\
& = \int_0^T \left|\left| \E \left[ \int_0^t x^2(s) \ d [M](s) \right] \right|\right|_Y \ d t \\
& \leq \int_0^T \left|\left| \E \left[ \int_0^T x^2(s) \ d [M](s) \right] \right|\right|_Y \ d t \\
& = T ||x||_{\mathcal M_M}^2.
\end{align*}
This shows that $\left(t \mapsto \int_0^t x(s) \ d M(s) \right) \in L^2([0,T];\mathcal M_T)$. Adaptedness follows from the definition of the It\^o integral.
\hfill $\square$
\end{quote}

\subsection{Lemma}
If $M=W$ is a Brownian motion, then $L^{2,a}([0,T];\mathcal M_T) \subset \mathcal M_W$.
\begin{quote}
\textsc{Proof:}
If $x \in L^2([0,T];\mathcal M_T)$, then
\begin{align*}
||x||_{\mathcal M_W}^2 & = \left|\left| \E \left[ \int_0^T x^2(s) \ d s \right] \right|\right|_Y = \left|\left| \int_0^T \E \left[ x^2(s) \right] \ d s \right|\right|_Y \\
& \leq \int_0^T \left|\left| \E \left[ x^2(s) \right] \right|\right|_Y \ d s = ||x||_{L^2([0,T];\mathcal M_T)}^2.
\end{align*}
\hfill $\square$
\end{quote}

\subsection{Corollary}
If $M=W$ is a Brownian motion and $x \in \mathcal M_W$, then the mapping $t \mapsto \int_0^t x(s) \ d W(s)$ is again an element of $\mathcal M_W$.
\begin{quote}
\textsc{Proof:}
This is an immediate consequence of the two previous results.
\hfill $\square$
\end{quote}

\subsection{Lemma}
\label{lem:functionalcommuteswithintegral1}
If $x \in L^{2,a}_M(\Omega \times [0,T]; X)$, and $\phi \in X^*$, then
\[ \phi\left(\int_0^T x(s) \ d M(s) \right) = \int_0^T \phi(x(s)) \ d M(s), \quad \textrm{a.s.}.\]

\begin{quote}
The proof is straightforward and omitted.
\end{quote}

\subsection{Lemma}
Let $x \in \mathcal M_M$, and $\phi \in X^*$. Suppose there exists a $\psi \in Y^*$ such that 
\[ \phi(\xi)^2 \leq \psi(\xi^2), \quad \textrm{for all} \ \xi \in X.\]
Then
\[ \phi \left( \int_0^T x(s) \ d M(s) \right) = \int_0^T \phi(x(s)) \ d M(s), \quad \textrm{a.s.}.\]
\begin{quote}
\textsc{Proof:} For $x \in \elem_M$ the claim is true. Suppose $(x_n)_{n \in \ints}$ is an approximating sequence of elementary functions for $x$ in $\mathcal M_M$.
The result follows from the estimate
\begin{align*}
& \E \left[ \left(\phi\left(\int_0^T x_n(s) \ d M(s) - \int_0^T x_m(s) \ d M(s)\right) \right)^2 \right] \\
& \leq \E \left[ \psi\left( \left(\int_0^T x_n(s) - x_m(s) \ d M(s)\right)^2 \right) \right] \\
& = \psi \left( \E \left[ \left(\int_0^T x_n(s) - x_m(s) \ d M(s)\right)^2 \right] \right) \\
& = \psi \left( \E \left[ \int_0^T (x_n-x_m)(s)^2 \ d[M](s) \right] \right) \\
& \leq ||\psi||_{Y^*} ||x_n - x_m||_{\mathcal M_M}^2.
\end{align*}
\hfill $\square$
\end{quote}

\subsection{Proposition}
\label{prop:stochintegralctsintoMt}Assume $M=W$ is a Brownian motion. Let $x \in \mathcal M_W$. The mapping
\[ t \mapsto \int_0^t x(s) \ d W(s) \]
is a continuous mapping from $[0,T]$ into $\mathcal M_T$.
\begin{quote}
\textsc{Proof:}
It suffices to show continuity at $t=0$.

First we will show the result for $x \in \elem_M$. Write
\[ x = \sum_{i=0}^n x_i \1_{[t_i, t_{i+1})}.\] For $t < t_1$, then
\[ \int_0^t x(s) \ d W(s) = x_0 W(t),\]
and, as $t \rightarrow 0$,
\[ \left|\left| \int_0^t x(s) \ d W(s) \right|\right|_{\mathcal M_T}^2 = \left|\left| \E \left[ x_0^2 W(t)^2 \right] \right|\right|_Y 
= \left|\left| \E \left[t x_0^2 \right] \right|\right|_Y \leq t ||x_0||_{L^2(\Omega,X)}^2 \rightarrow 0. \]

Now for $x \in \mathcal M_W$, choose $y \in \elem_W$ such that
\[ ||x - y||_{\mathcal M_W} < \epsilon /2.\]

Choose $t > 0$ such that
\[ \left|\left|\int_0^t y(s) \ d W(s) \right|\right|_{\mathcal M_T} < \epsilon / 2.\]

Then
\begin{align*} \left|\left| \int_0^t x(s) - y(s) \ d W(s) \right|\right|_{\mathcal M_T}^2 & = \left|\left| \E \left[ \left( \int_0^t x(s) - y(s) \ d W(s) \right)^2 \right] \right|\right|_Y \\
& = \left|\left| \E \left[ \int_0^t (x(s) - y(s))^2 \ d s \right] \right|\right|_Y \\
& \leq \left|\left| \E \left[ \int_0^T (x(s) - y(s))^2 \ d s \right] \right|\right|_Y = || x - y ||_{\mathcal M_W}^2 < (\epsilon / 2)^2,
\end{align*}
so
\[ \left|\left| \int_0^t x(s) d W(s) \right|\right|_{\mathcal M_T} \leq \left|\left| \int_0^t x(s) - y(s) \ d W(s) \right|\right|_{\mathcal M_T} + \left|\left| \int_0^t y(s) \ d W(s) \right|\right|_{\mathcal M_T} < \epsilon.\]
\hfill $\square$
\end{quote}

\section{The It\^o integral in $C(K)$}
Throughout this section, let $(\Omega, \mathcal F, \P)$ be a probability space, $(\mathcal F_t)_{t \geq 0}$ a filtration on $\Omega$, $K$ a compact metric space, and $X = Y = C(K)$. We say a process is adapted when it is adapted to $(\mathcal F_t)_{t \geq 0}$. 

\subsection{Proposition} \label{prop:characterizationCK} For $\mathcal M_M$ as defined in Definition~\ref{def:stochintegrable}, with $X, Y = C(K)$, 
\[ \mathcal M_M \cong C(K; L^{2,a}_M(\Omega \times [0,T])).\]

\subsubsection{Lemma} \label{lem:MnormCK} Under the assumptions of Proposition~\ref{prop:characterizationCK}, let $x \in L^{2,a}_M(\Omega \times [0,T];C(K))$. Then 
\begin{align*} || x || & := \left(\sup_{k \in K} \E \left[ \int_0^T x(t,\cdot,k)^2 \ d[M](t) \right] \right)^{\half} \\
& = \sup_{k \in K} ||x (\cdot,\cdot,k) ||_{L^2_M(\Omega \times([0,T])}.
\end{align*} \hfill $\square$

\subsubsection{Lemma} After a change of variable order,
\label{lem:L2inCK}
\[ L^{2,a}_M(\Omega \times [0,T];C(K)) \subset C(K;L^{2,a}_M(\Omega \times [0,T])).\]
\begin{quote} \textsc{Proof:}
Let $x \in L^{2,a}_M(\Omega \times [0,T];C(K))$. By the previous lemma, and by~(\ref{eq:estimateMnorm}),
\[ ||x||_{C(K;L^2_M(\Omega \times [0,T]))} = \sup_{k \in K} ||x(\cdot, \cdot, k)||_{L^2_M(\Omega \times [0,T])} = ||x|| \leq ||x||_{L^2_M(\Omega \times [0,T];C(K))} < \infty.\]

Let $k \in K$. For almost all $(\omega,t) \in \Omega \times [0,T]$, $|x(t,\omega,k) - x(t,\omega,l)| \rightarrow 0$ as $l \rightarrow k$ in $K$. Because $|x(t,\omega,k) - x(t,\omega,l)| \leq 2 ||x(t,\omega,\cdot)||_K$, and $\int_0^T \E [ ||x(t)||_K^2 ] \ d t < \infty$, by dominated convergence 
\[ \E \left[ \int_0^T ||x(t,\cdot,k) - x(t,\cdot,l)||_K^2 \ d[M](t) \right] \rightarrow 0, \quad l \rightarrow k,\]
so $x \in C(K; L^2_M(\Omega \times [0,T]))$.

For almost all $t \in [0,T]$ we have that $x(t) \in L^2(\Omega; C(K))$ is $\mathcal F_t$-measurable. Because measurability implies weak measurability, $x(t,\cdot,k)$ is $\mathcal F_t$-measurable for all $k \in K$. This implies that, for all $k \in K$, $x(t,\cdot,k)$ is $\mathcal F_t$-measurable for almost all $t \in [0,T]$, i.e. for all $k \in K$ there exists an $(\mathcal F_t)$-adapted version of $x(\cdot,\cdot,k)$.
\hfill $\square$
\end{quote}

\subsubsection{Lemma} 
\label{lem:cauchyinCK}
The space $C(K;L^{2,a}_M(\Omega \times [0,T]))$ is complete.

\begin{quote}
\textsc{Proof:}
Let $(x_n)_{n \in \ints}$ be a sequence in $C(K;L^{2,a}_M(\Omega \times [0,T]))$ .
Because $C(K;L^2_M(\Omega \times [0,T]))$ is complete, $(x_n)_{n \in \ints}$ has a limit $x \in C(K;L^2_M(\Omega \times [0,T]))$.

For $k \in K$, $x_n(k) \rightarrow x(k)$ in $L^2_M(\Omega \times [0,T])$. This implies $\Omega \times [0,T]$-almost everywhere convergence of a subsequence, $(x_{n_l}(k))$ say. For almost all $t \in [0,T]$, $x_{n_l}(k,t) \in L^2(\Omega)$ is $\mathcal F_t$-measurable, and so is its limit $x(k,t)$. This shows that $x(k)$ is $(\mathcal F_t)$-adapted for all $k \in K$.
\hfill $\square$
\end{quote}

\subsubsection{Lemma}
\label{lem:L2SCKdenseinCKL2S}
Let $(S, \Sigma, \mu)$ be a measure space and let $K$ be a compact set. Then
\[ \{ y \in L^2(S;C(K)): y = \sum_{i=1}^n \alpha_i y_i, \textrm{with} \ \alpha_i \in C(K), y_i \in L^2(S), i = 1, \hdots, n \} \]
is dense in $C(K;L^2(S))$.
\begin{quote}
\textsc{Proof:} 
If $y \in L^2(S;C(K))$, then
\[ ||y||_{C(K;L^2(S))}^2 = \sup_{k \in K} || y(k) ||_{L^2(S)}^2 \leq \int_S \sup_{k \in K} y(k)^2 \ d \mu = ||y||_{L^2(S;C(K))}^2 < \infty.\]

Let $x \in C(K;L^2(S))$. Choose an $\epsilon > 0$.
Because $K$ is compact, $x$ takes it values in a compact subset of $L^2(S)$. Since the collection of balls 
\[ B(y,\epsilon/2) = \{ \xi \in L^2(S) : ||\xi - y||_{L^2(S)} < \epsilon /2 \}, \quad y \in x(K), \]
is an open cover of $x(K)$, there exists a finite set $y_1, \hdots, y_n \in x(K)$ (where we may assyme that $y_i \neq 0$, $i = 1, \hdots, n$), such that
\begin{equation} \label{eq:xKcover} x(K) \subset \bigcup_{i=1}^n B(y_i,\epsilon/3).\end{equation}
Now $\mathrm{span}(y_1, \hdots, y_n)$ is a closed linear subspace of $x(K)$. Let $P : L^2(S) \rightarrow \mathrm{span}(y_1,\hdots, y_n)$ denote the orthogonal projection from $L^2(S)$ onto $\mathrm{span}(y_1,\hdots,y_n)$. Then for $x \in L^2(S)$,
\[ P x = \sum_{i=1}^n \frac{\left< x, y_i \right>}{ || y_i||^2_{L^2(S)} } y_i \]
where $\left<\cdot,\cdot\right>$ denotes the inner product on $L^2(S)$. Define
\[ y : K \mapsto L^2(S), \quad k \mapsto P [x(k)] = \sum_{i=1}^n  \alpha_i(k) y_i,\]
with 
\[ \alpha_i : k \mapsto \frac{\left< x(k), y_i \right>}{ || y_i||^2_{L^2(S)}}, \quad i = 1, \hdots, n.\]

Furthermore, for $k \in K$, let $j : K \rightarrow \{1, \hdots, n\}$ such that $x(k) \in B(y_{j(k)}, \epsilon /2)$ (possible because of~(\ref{eq:xKcover})). Then
\begin{align*}
|| x(k) - y(k) ||_{L^2(S)} & \leq ||x(k) - y_{j(k)} ||_{L^2(S)} + ||y_{j(k)} - y(k) ||_{L^2(S)} \\
& = ||x(k) - y_{j(k)} ||_{L^2(S)} + ||P [y_{j(k)}] - P [x(k)] ||_{L^2(S)} \\
& \leq 2 ||x(k) - y_{j(k)} ||_{L^2(S)} < \epsilon.
\end{align*}
So 
\[ \sup_{k \in K} ||x(k) - y(k) ||_{L^2(S)} < \epsilon.\]

For $i = 1,\hdots, n$ and $k, l \in K$,
\begin{align*}
| \alpha_i(k) - \alpha_i(l)| & = \left| \sum_{i=1}^n  \frac{\left< x(k) - x(l), y_i \right>}{ || y_i||^2_{L^2(S)}} \right| \\
& \leq \sum_{i=1}^n \frac{|| x(k) - x(l)||_{L^2(S)}}{||y_i||_{L^2(S)} }.
\end{align*}
So, because $x \in C(K;L^2(S))$, we see that $\alpha_i \in C(K), i=1, \hdots, n$. So
\[ y : S \rightarrow C(K), \quad m \rightarrow \sum_{i=1}^n y_i(m) \alpha_i,\]
has the correct form, is strongly measurable, and 
\begin{align*} || y||_{L^2(S;C(K))}^2 & = \int_S \left|\left|\sum_{i=1}^n \alpha_i y(m) \right|\right|_{C(K)}^2 \ d \mu(m) \\
& \leq \sum_{i=1}^n ||\alpha_i||_{C(K)}^2 \int_S |y_i|^2 \ d \mu = \sum_{i=1}^n ||\alpha_i||_{C(K)}^2 ||y_i||_{L^2(S)}^2 < \infty.
\end{align*}
\hfill $\square$
\end{quote}

\subsubsection{Lemma}
\label{lem:elemdenseinCK}
$\elem([0,T]; L^2(\Omega; C(K)))$ is dense in 
\[ \{ x \in C(K;L^2([0,T] \times \Omega)) : x \ \textrm{is adapted to} \ (\mathcal F_t) \ \textrm{for all} \ k \in K \}.\]

\begin{quote}
\textsc{Proof:}
Let $x \in C(K; L^2([0,T]\times\Omega))$ be adapted to $(\mathcal F_t)$ for all $k \in K$.
By Lemma~\ref{lem:L2SCKdenseinCKL2S}, there is an $y \in L^2([0,T] \times \Omega ; C(K) )$ such that
\[ y = \sum_{i=1}^n \alpha_i y_i,\] with $\alpha_i \in C(K)$ and $y_i \in L^2([0,T]\times\Omega)$ for $i=1,\hdots,n$, and such that
\[ || y - x ||_{C(K;L^2([0,T]\times\Omega))} < \epsilon / 2.\]
Because the adapted functions form a closed subspace of $L^2([0,T]\times\Omega)$, there exists an orthogonal projection operator
\[ Q : L^2([0,T]\times \Omega) \rightarrow \{ x \in L^2([0,T]\times\Omega) : \ x \ \textrm{adapted to} \ (\mathcal F_t)_{t \in [0,T]} \}.\]
When we apply Proposition~\ref{lem:elemdenseinL2} to $X = \reals$, this shows that $\elem([0,T];L^2(\Omega))$ is dense in $\{ x \in L^2([0,T]\times\Omega): x \ \textrm{adapted to} \ (\mathcal F_t)_{t\in[0,T]} \}$. Therefore, for $i = 1,\hdots, n$, there exists an $z_i \in \elem([0,T];L^2(\Omega))$ such that 
\[ || z_i - Q y_i ||_{L^2([0,T]\times\Omega)} < \frac {\epsilon}{2 n (||\alpha_i||_{C(K)} + 1)}.\]

Then, for
\[ z = \sum_{i=1}^n \alpha_i z_i,\]
we find that
\[ || z - x ||_{C(K;L^2([0,T]\times\Omega))} \leq ||z - \sum_{i=1}^n \alpha_i Q y_i ||_{C(K;L^2([0,T] \times \Omega))} + || \sum_{i=1}^n \alpha_i Q y_i - x ||_{C(K;L^2([0,T]\times\Omega))},\]
with
\begin{align*} || z - \sum_{i=1}^n \alpha_i Q y_i ||_{C(K;L^2([0,T]\times\Omega))} & = || \sum_{i=1}^n \alpha_i (z_i - Q y_i) ||_{C(K;L^2([0,T]\times\Omega))} \\
& \leq \sup_{k \in K} \sum_{i=1}^n |\alpha_i(k)| \ || z_i - Q y_i ||_{L^2([0,T]\times\Omega)} \\
& \leq \sum_{i=1}^n ||\alpha_i||_{C(K)} || z_i - Q y_i ||_{L^2([0,T]\times\Omega)} < \epsilon / 2,
\end{align*}
and
\begin{align*}
|| \sum_{i=1}^n \alpha_i Q y_i - x ||_{C(K;L^2([0,T]\times\Omega)} & = \sup_{k \in K} || Q [y(k)] - Q [x(k)] ||_{L^2([0,T]\times\Omega)} \\
& \leq \sup_{k \in K} || y(k) - x(k)||_{L^2([0,T]\times\Omega)} < \epsilon / 2.
\end{align*}
So $|| x - z ||_{C(K;L^2([0,T]\times\Omega))} < \epsilon$, and $z \in \elem([0,T]; L^2(\Omega; C(K))$.
\hfill $\square$
\end{quote}

\begin{quote}
\textsc{Proof of Proposition~\ref{prop:characterizationCK} }: The desired result is obtained from Lemma~\ref{lem:L2inCK}, Lemma~\ref{lem:cauchyinCK} and Lemma~\ref{lem:elemdenseinCK}.
\end{quote}

\subsection{Corollary}
Under the assumptions of Proposition~\ref{prop:characterizationCK}, for $t \in [0,T]$
\[ R_t \cong C(K; L^2((\Omega, \mathcal F_t))) .\]

\begin{quote}
\textsc{Proof:} The proof of this statement is a simplified version of the proof of Proposition~\ref{prop:characterizationCK}. \hfill $\square$
\end{quote}

\subsection{Proposition}
Let $K = [a,b]$, a closed interval in $\reals$, and let $X = C([a,b])$. Suppose $x \in R \subset C([a,b], L^2([0,T] \times \Omega))$, and $t \in [0,T$], and
\[ |x(\tau)(s,\omega) - x(\sigma)(s,\omega)| \leq L(\omega) |\tau - \sigma|^\beta,\]
with $\E [L^2] < \infty$ and $\beta > \half$, for all $s \in [0,t]$, $\omega \in \Omega$.

Then there exists a version of
\[ \sigma \mapsto \int_0^t x(s) \ d W(s)(\sigma) \] 
that is H\"older continuous on $[a,b]$ with exponent $\gamma$, for every $\gamma \in (0, \beta - \half)$.

\begin{quote}
\textsc{Proof:} By the mentioned assumptions,
\begin{align*}
\E \left[ \left( \int_0^t x(s) \ d W(s) (\tau) - \int_0^t x(s) \ d W(s) (\sigma) \right)^2 \right] & = \E \left[ \int_0^t (x(s)(\tau) - x(s)(\sigma))^2 \ d s \right] \\
& \leq \E \left[ L^2 \int_0^t |\tau - \sigma|^{2 \beta} \ d s \right] = t \E[L^2]\ |\tau - \sigma|^{2 \beta}.
\end{align*}

Then Kolmogorov's continuity criterion (\cite{karatzasshreve}, Theorem 2.2.8), states that 
\[ \sigma \mapsto \int_0^t x(s)(\sigma) \ d W(s) \] has a H\"older continuous version with exponent $\gamma$, for every $\gamma \in (0, \beta - \half)$.
\hfill $\square$
\end{quote}

\section{The It\^o integral in $L^p(\mu)$}
\subsection{Proposition}
\label{prop:characterizationLp}
Let $(S,\Sigma,\mu)$ be a $\sigma$-finite measure space. Let $2 \leq p \leq \infty$. Suppose $X = L^p(S;\reals)$ and $Y = L^{p/2}(\mu)$ (possibly with $p = p/2 = \infty$). Then
\[ R \cong \{ x \in L^p((S;\mu), L^2([0,T] \times \Omega) ) : x \ \textrm{adapted to} \  (\mathcal F_t), \ \mu \textrm{-almost everywhere} \} .\]

To prove this, we need a few lemmas.

\subsubsection{Lemma} In case $X = L^p(\mu)$ and $Y = L^{p/2}(\mu)$ for $p \geq 2$, then for $x \in \elem([0,T];L^2(\Omega;X))$ we have
\begin{align*} || x ||_R & = \left|\left| \int_0^T \E [ x^2(t)] \ d t \right|\right|_Y^{\half} = \left( \int_S \left( \int_0^T \E [x^2(t)] \ d t \right)^{p/2} \ d \mu \right)^{1/p} \\
& = \left( \int_S || x ||_{L^2([0,T] \times \Omega)}^p \ d \mu \right)^{1/p} = ||x ||_{L^p(S; L^2([0,T] \times \Omega))}.
\end{align*} \hfill $\square$

\subsubsection{Lemma} 
\label{lem:eleminLp}
After a change of variable order,
\[ \elem([0,T];L^2(\Omega; L^p(S))) \subset \left\{ x \in L^p(S;L^2([0,T] \times \Omega)) : x \ \textrm{adapted to} \ (\mathcal F_t), \mu \textrm{-almost everywhere} \right\} .\]

\begin{quote} \textsc{Proof:}
Let $x \in \elem([0,T]; L^2(\Omega; L^p(S)))$. Then by the previous lemma, $||x||_R = ||x||_{L^p(S;L^2([0,T] \times \Omega))} < \infty$, so $x \in L^p(S;L^2([0,T] \times \Omega))$. Suppose first that $\mu(S) < \infty$.

This $x$ can be written as
\[ x = \sum_{i=0}^n x_i \1_{[t_i, t_{i+1})} \1_{S_m}, \]
where all the $x_i \in L^2(\Omega; L^p(S))$ are $\mathcal F_{t_i}$-measurable.

Because $x_i \in L^2(\Omega; L^p(S)) \subset L^1(\Omega; L^1( S)) \cong L^1(\Omega \times S)$, we see that $\omega \mapsto x_i(\omega,s)$ is $\mathcal F_{t_i}$-measurable for $\mu$-almost all $s \in S$. So $(s, t, \omega) \mapsto x_i(\omega, s) \1_{[t_i, t_{i+1})}(t)$ is $(\mathcal F_t)$-adapted, for $\mu$-almost all $s$, for all $i = 0, \hdots, n$, and so is $(s, t, \omega) \mapsto x(t,\omega, s)$.

For $\mu(S) = \infty$, let $S = \cup_{m=1}^\infty S_m$ with $\mu(S_m) < \infty$ for all $m$. A similar argument shows that $x_i \1_{S_m}$ is $\mathcal F_{t_i}$-measurable for $\mu$-almost all $s$ and all $m \in \ints$, so $x_i$ is $\mathcal F_{t_i}$-measurable for $\mu$-almost all $s$, and the same result follows.
\hfill $\square$
\end{quote}

\subsubsection{Lemma} 
\label{lem:cauchyinLp}
Let $(x_n)$ be a Cauchy sequence of $\mu$-almost everywhere $(\mathcal F_t)$-adapted functions in $L^p(S;L^2([0,T] \times \Omega))$. Then its limit $x$ exists and is $(\mathcal F_t)$-adapted.

\begin{quote} \textsc{Proof:} Since $L^p(S; L^2([0,T] \times \Omega))$ is a Banach space, the limit $x = \lim_{n \rightarrow \infty} x_n$ exists.

Convergence in $L^p$ implies $\mu$-almost everywhere convergence of a subsequence $(x_{n_k})$. So $x_{n_k}(s) \rightarrow x(s)$ in $L^2([0,T] \times \Omega)$ for $\mu$-almost all $s \in S$. This implies the existence of a further subsequence $(x_{n_{k_l}})$ so that $x_{n_{k_l}}(s,t) \rightarrow x(s,t)$ for almost all $t \in [0,T]$, and $\mu$-almost all $s \in S$. Since $x_{n_{k_l}}(s)$ is $(\mathcal F_t)$-adapted for $\mu$-almost all $s$, $x_{n_{k_l}}(s,t)$ is $\mathcal F_t$-measurable for almost all $t \in [0,T]$, and $\mu$-almost all $s \in S$, and so is the limit $x(s,t)$. This shows that $x$ is $(\mathcal F_t)$-adapted, $\mu$-almost everywhere. \hfill $\square$
\end{quote}

\subsubsection{Lemma} 
\label{lem:elemdenseinLp}
$\elem([0,T]; L^2(\Omega; L^p(S)))$ is dense in 
\[ \{ x \in L^p(S;L^2([0,T] \times \Omega)) : x \ \textrm{is adapted to} \ (\mathcal F_t), \mu\textrm{-almost everywhere} \}.\]

\begin{quote}
\textsc{Proof:} Let $x \in L^p(S;L^2([0,T] \times \Omega))$, adapted to $(\mathcal F_t)$ and choose $\epsilon > 0$. There exists a simple function $y \in \simple(S,L^2([0,T] \times \Omega))$ such that
\[ || x - y||_{L^p(S; L^2([0,T] \times \Omega))} < \epsilon / 2.\]
Write \[ y = \sum_{i = 1}^n y_i \1_{A_i},\] with, for all $i = 1, \hdots, n$,
$y_i \in L^2([0,T] \times \Omega)$ and $A_i \subseteq S_0$, with $\mu(S_0) < \infty$ (possible because of $\sigma$-finiteness of $S$).
The set of $(\mathcal F_t)$-adapted functions in $L^2([0,T] \times \Omega)$ is a closed subspace of $L^2([0,T] \times \Omega)$, so we can define the orthogonal projection 
\[ P : L^2([0,T] \times \Omega) \rightarrow \{ x \in L^2([0,T] \times \Omega): x \ \textrm{adapted to} \ (\mathcal F_t) \}. \]
Define \[ z : = \sum_{i=1}^n P y_i \1_{A_i}. \]
Then
\begin{align*} ||x - z||^p_{L^p(S;L^2([0,T] \times \Omega))} & = \int_S ||x(s) - z(s) ||_{L^2([0,T] \times \Omega)}^p \ d \mu(s) \\
& = \int_S || P (x(s) - y(s)) ||_{L^2([0,T] \times \Omega)}^p \ d \mu(s) \\
& \leq \int_S || x(s) - y(s) ||_{L^2([0,T] \times \Omega)}^p \ d \mu(s) < \left( \frac \epsilon 2 \right)^p.
\end{align*}
Because $P y_i \in L^2([0,T] \times \Omega)$ is $(\mathcal F_t)$-adapted for $i = 1, \hdots, n$, by Lemma~\ref{lem:elemdenseinL2} there exists an elementary function $a_i \in \elem([0,T]; L^2(\Omega))$ such that
\[ || P y_i - a_i ||_{L^2([0,T] \times \Omega)} < \frac \epsilon 2 \left( \frac 1 {n \mu(S_0) + 1} \right)^{1/p}. \]
Then
\[ \sum_{i=1}^n a_i \1_{A_i} \in \elem([0,T]; L^2(\Omega; L^p(S))),\] and
\begin{align*} | | \sum_{i=1}^n a_i \1_{A_i} - z  ||_{L^p(S;L^2([0,T]\times \Omega))}^p & = \int_S \sum_{i=1}^n || a_i - P y_i ||_{L^2([0,T]\times\Omega)}^p \1_{A_i} \ d \mu \\
& < \left( \frac \epsilon 3 \right)^p \sum_{i=1}^n \frac 1 {n \mu(S_0) + 1} \mu(A_i) < \left( \frac \epsilon 2 \right)^p.
\end{align*}
Hence 
\begin{align*} || \sum_{i=1}^n a_i \1_{A_i} - x ||_{L^p(S;L^2([0,T] \times \Omega))} & \leq || \sum_{i=1}^n a_i \1_{A_i} - z ||_{L^p(\hdots)} + || z - x ||_{L^p(\hdots)} < \frac \epsilon 2 + \frac \epsilon 2 = \epsilon.
\end{align*}
\hfill $\square$
\end{quote}

\begin{quote}
\textsc{Proof of Proposition~\ref{prop:characterizationLp}:} The result follows from Lemmas~\ref{lem:eleminLp},~\ref{lem:cauchyinLp} and~\ref{lem:elemdenseinLp}.
\hfill $\square$
\end{quote}

\subsubsection*{Remark} Note the similarities between Propositions~\ref{prop:characterizationCK} and ~\ref{prop:characterizationLp}, Lemmas~\ref{lem:L2inCK} and~\ref{lem:eleminLp}, Lemmas~\ref{lem:cauchyinCK} and ~\ref{lem:cauchyinLp}, and Lemmas~\ref{lem:elemdenseinCK} and ~\ref{lem:elemdenseinLp}.

\subsection{Lemma}
\label{lem:evaluationcommutes}
Let $X$ be a function space. Then for $t \in [0,T]$ and $x \in R$,
\[ \int_0^t x(s)(\cdot) \ d W(s) = \int_0^t x(s) \ d W(s) (\cdot), \quad \textrm{and} \quad \int_0^t x(s)(\cdot) \ d s = \int_0^t x(s) \ d s (\cdot).\]

\begin{quote}
\textsc{Proof:}
Let $(x_n)_{n \in \ints}$ be an approximation of $x$ by elementary functions, $x_n \in \elem([0,T]; L^2(\Omega;X))$, $n \in \ints$. Then for all $n \in \ints$,

\[ \int_0^t x_n(s)(\cdot) \ d W(s) = \int_0^t x(s) \ d W(s) (\cdot) ,\] because here the stochastic integral is just addition of a finite number of terms.

By the classical, real-valued It\^o isometry,
\begin{align*}
\left|\left| \int_0^t x_n(s) (\cdot) \ d W(s) - \int_0^t x_m(s) (\cdot) \ d W(s) \right|\right|_{R_t}^2 & = \left|\left| \E \left[ \left(\int_0^t x_n(s) (\cdot) - x_m(s) (\cdot) \ d W(s) \right)^2 \right] \right|\right|_Y \\
& = \left|\left| \E \left[ \int_0^t (x_n(s) -x_m(s))^2(\cdot) \ d s \right] \right|\right|_Y \\
& = \left|\left| \E \left[ \int_0^t (x_n(s) - x_m(s) )^2 \ d s \right] (\cdot) \right|\right|_Y \\
& \leq ||x_n - x_m||_{R}^2,
\end{align*}
which shows that as $n \rightarrow \infty$, the sequence
\[ \int_0^t x_n(s) (\cdot) \ d W(s) \] converges to a limit in $R_t$, necessarily equal to \[ \int_0^t x(s) \ d W(s)(\cdot) .\]

The proof for point evaluation of the integral $\int_0^t x(s) \ d s$ is similar.
\hfill $\square$
\end{quote}

\subsection{Proposition}
If $X = L^p(S)$, and $x \in R$, then for almost all $\sigma \in S$, $\omega \in \Omega$, the mapping
\[ [0,T] \rightarrow \reals, \quad t \mapsto  \left(\int_0^t x(s) \ d W(s) \right)(\sigma, \omega) \]
is continuous.
\begin{quote}
\textsc{Proof:}
Since $x \in R \cong L^p(S, L^2([0,T] \times \Omega))$, we have that $||x(\sigma)||_{L^2([0,T] \times \Omega)} < \infty$ for almost all $\sigma \in S$.
Hence
\[ \P \left( \int_0^T x^2(s) (\sigma) \ d s < \infty \right) < 1, \quad \textrm{a.a.} \ \sigma \in S.\]
Therefore (by \cite{karatzasshreve}, Definition 3.2.23) the real-valued stochastic integral,
\[ t \mapsto \left(\int_0^t x(s)(\sigma) d W(s)\right)(\omega),\] is continuous for almost all $\sigma \in S$ and $\omega \in \Omega$, and by Lemma~\ref{lem:evaluationcommutes} this is equivalent to saying that
\[ t \mapsto \left(\int_0^t x(s) d W(s)\right) (\sigma, \omega) \] is continuous for almost all $\sigma \in S$, $\omega \in \Omega$.
\hfill $\square$
\end{quote}

\appendix

\section{Vector-valued integration}
\label{app:integration}

Throughout this section, let $(S,\Sigma,\mu)$ denote a $\sigma$-finite measure space, and $X$ a real Banach space with dual space $X^*$. 

A function $f: S \rightarrow X$ is called \emph{measurable} if $f^\leftarrow (V) \in \Sigma$ for every Borel set $V$ in $X$.

A function $f: S \rightarrow X$ is called \emph{weakly measurable} if $\phi(f)$ is measurable for every $\phi \in X^*$.

\subsection{Weak integral}
A function $f: S \rightarrow X$ is called \emph{weakly integrable} if 
\begin{itemize}
\item[(i)] $f$ is weakly measurable,
\item[(ii)] for every $\phi \in X^*$ the integral
\[ \int_S \phi(f) \ d \mu\] exists, and
\item[(iii)] there exists an $I(f) \in X$ such that
\[ \int_S \phi(f) \ d \mu = \phi(I(f)) \quad \textrm{for all} \ \phi \in X^*. \] 
\end{itemize}
In this case $I(f)$ is called the \emph{weak integral} of $f$.

Furthermore we write
\[ \int_S f d \mu := I(f). \]

By Hahn-Banach the weak integral is uniquely defined.

\subsection{Strong integral}
A function $f: S \rightarrow X$ is called \emph{simple} if there exist $E_1, \hdots, E_n \in \Sigma$ with $\P(E_k) < \infty$ for $1 \leq k \leq n$ and $f_1, \hdots f_n \in X$ such that 
\[ f(s) = \sum_{k=1}^n \1_{E_k}(s) f_k \quad \textrm{for all} \ s \in S.\]

A function $f: S \rightarrow X$ is called \emph{strongly measurable} if there exists a sequence of simple functions $(f_n)$ such that $f_n(x) \rightarrow f(x)$ for all $x \in S$.

A function $f: S \rightarrow X$ is called \emph{$\mu$-strongly measurable} if it is $\mu$-almost everywhere equal to a strongly measurable function.

A $\mu$-strongly measurable function $f: S \rightarrow X$ is called \emph{strongly integrable} if there exists a sequence $(f_n)$ of simple functions such that
\[ \int_S ||f - f_n|| \ d \mu \rightarrow 0, \quad n \rightarrow \infty.\]
Then $(\int_S f_n \ d \mu)$ is a Cauchy sequence in $X$ and we can define the \emph{strong integral} of $f$ by
\[ \int_S f \ d \mu := \lim_{n \rightarrow \infty} \int_S f_n \ d \mu.\]

Some results on weak and strong integrals are listed below.

\begin{itemize}
\item[(i)] A $\mu$-strongly measurable function $f: S \rightarrow X$ is strongly integrable if and only if $\int_S ||f|| \ d \mu < \infty$.
\item[(ii)] If $X$ is separable, then a function $f: S \rightarrow X$ is strongly integrable if and only if it is weakly measurable and $\int_S ||f|| \ d \mu < \infty$.
\item[(iii)] A function $f : S \rightarrow X$ is strongly measurable if and only if $f$ is separably valued and weakly measurable.
\end{itemize}

\subsection{$L^p$ spaces}
For $1 \leq p < \infty$ let $\mathcal L^p(\mu;X)$ denote the space of all strongly integrable functions $f: S \rightarrow X$ with $\int_S ||f||^p < \infty$.
For $f \in \mathcal L^p(\mu;X)$ define
\[ ||f||_p := ||f||_{L^p(\mu;X)} = \left( \int_S ||f||^p d \mu \right)^{\frac 1 p}.\]
Then $||\cdot||_p$ is a seminorm on $\mathcal L^p(\mu;X)$ whose kernel consists of the functions that are almost everywhere equal to 0. Let $L^p(\mu;X)$ denote the quotient space of $\mathcal L^1(\mu;X)$ over the kernel of $||\cdot||_p$ and denote its norm again by $||\cdot||_p = ||\cdot||_{L^p(\mu;X)}$.

The space $L^\infty(\mu,X)$ consists of all equivalence classes of $\mu$-strongly measurable $f: S \rightarrow X$ such that $||f(\cdot)||_X$ is bounded on a $\mu$-full set, endowed with the norm $||f||_\infty = || ||f(\cdot)||_X ||_\infty, f \in L^\infty(\mu,X)$.

The spaces $(L^p(\mu;X), ||\cdot||_p)$, $1 \leq p \leq \infty$ are Banach spaces.

As for the real-valued $L^p$-spaces, we will use the notations 
\[ L^p(\mu; X) = L^p( (S, \mu); X) = L^p( S; X), \] dependent on the context in which they occur.

\subsection{Proposition}
\label{prop:bochner=pointwise}
If $X$ is a function space, or a space of equivalence classes of functions, and $f \in L^1(\mu; X)$ then
\[ \left(\int_S f \ d \mu \right)(\cdot) = \int_S f(\cdot) \ d \mu \quad \textrm{in} \ X, \]
i.e. evaluating the Bochner integral is equivalent to pointwise integration.
\begin{quote}
\textsc{Proof:}
Let $(f_n)$ be a sequence of simple functions so that $\int_S ||f_n - f||_X \ d \mu \rightarrow 0$ as $n \rightarrow \infty$.
For all $f_n, n \in \ints$,
\[ \left(\int_S f_n \ d \mu\right)(\cdot) = \int_S f_n(\cdot) \ d \mu,\] because here integration is just a summation.
Then
\begin{align*}
|| \left(\int_S f \ d \mu \right) (\cdot) - \int_S f(\cdot) \ d \mu ||_X & = || \lim_{n \rightarrow \infty} \left(\int_S f_n \ d \mu \right)(\cdot) - \lim_{n \rightarrow \infty} \int_S f_n(\cdot) \ d \mu ||_X \\
& = \lim_{n \rightarrow \infty} || \left( \int_S f_n \ d \mu \right) (\cdot) - \int_S f_n(\cdot) \ d \mu ||_X = 0. \end{align*}
\hfill $\square$
\end{quote}

\subsection{Banach lattice}
\label{app:Banachlattice}
An \emph{ordered vector space} is a vector space $X$ over $\reals$ with a partial order $\preceq$ such that if $x,y \in X$ and $x \preceq y$, then
\begin{itemize}
\item[(i)] $x+z \preceq y+z$ whenever $z \in X$; and 
\item[(ii)] $tx \preceq ty$ whenever $t > 0$.
\end{itemize}

A \emph{vector lattice} is an ordered vector space $X$ such that
\begin{itemize}
\item[(iii)] every pair of elements $x,y$ in $X$ has a least upper bound $x \vee y$ in $X$.
\end{itemize}

For every element $x$ of a vector lattice, let the \emph{absolute value} of $x$ be defined by the formula $|x| = x \vee (-x)$. A real normed space $X$ that is also a vector lattice is a \emph{normed lattice} if
\begin{itemize}
\item[(iv)]
$||x|| \leq ||y||$ whenever $x,y \in X$ and $|x| \preceq |y|$.
\end{itemize}

A \emph{Banach lattice} is a real Banach space that is a normed lattice.

\subsection{Properties of vector lattices}
Let $a, b, c \in X$. 
\begin{itemize}
\item[(i)] If $a \succeq 0$, then $-a \preceq 0$;
\item[(ii)] $a + (b \vee c) = (a + b) \vee (a + c)$;
\item[(iii)] $|a| \geq 0$;
\item[(iv)] $a \succeq 0$ if and only if $|a| = a$;
\item[(v)] $|| \ |x| \ ||_X = ||x||_X$, so $|\cdot| : X \rightarrow X$ is continuous.
\end{itemize}

\subsection{Lemma}
\label{lem:integralpositive}
If $X$ is a Banach lattice, and $x \in L^1(M; X)$, with $x \succeq 0$ on $A \subset M$, then
\[ \int_M x \1_A \ d \mu \succeq 0.\]

\begin{quote}
\textsc{Proof:}
It suffices to consider functions $x \in L^1(M;X)$ such that $x \succeq 0$ on $M$, because if $x \succeq 0$ on $A \subset M$, then $x \1_A \succeq 0$ on $M$.

Let $(x_n)_{n \in \ints}$ be a sequence of simple functions such that $x_n \rightarrow x$, $\mu$-almost everywhere, and such that
\[ \int_M ||x_n - x||_X \ d\mu \rightarrow 0.\]
For all $n \in \ints$, define $\tilde x_n : = x_n \vee 0$. Then $\tilde x_n \succeq 0$ is also a simple function.

Using some elementary properties of partial ordering and absolute value,
\[ |\tilde x_n - x | \preceq |x_n - x|, \quad \mu\textrm{-almost everywhere},\] 
so
\[ ||\tilde x_n - x ||_X \leq ||x_n - x||_X, \quad \mu\textrm{-almost everywhere}.\]

For all $n \in \ints$
\[ \int_M \tilde x_n \ d \mu \geq 0, \]
because the lefthand side is just a summation of positive elements.

Because $|\cdot| : X \rightarrow X$ is continuous,
\begin{align*} \int_M x \1_A \ d \mu & = \lim_{n \rightarrow \infty} \int_M \tilde x_n \1_A \ d \mu = \lim_{n \rightarrow \infty} \left| \int_M \tilde x_n \1_A \ d \mu \right| \\
& = \left| \lim_{n \rightarrow \infty} \int_M \tilde x_n \1_A \ d \mu \right| = \left| \int_M x \1_A \ d \mu \right|,
\end{align*}
so that
\[ \int_M x \1_A \ d \mu \succeq 0.\]
\hfill $\square$
\end{quote}

\subsection{Lemma}
\label{lem:normingsequence}
Let $Y$ be a separable closed subspace of $X$. Then $X^*$ containts a sequence $(\phi_n)_{n \in \ints}$ with $||\phi_n||_{X^*} = 1$ for all $n \in \ints$, which is norming for $Y$, that is,
\[ ||y||_X = \sup_{n \in \ints} |\phi_n(y)| \quad \textrm{for all} \ y \in Y.\]

\begin{quote}
\textsc{Proof:} Choose a dense sequence $(x_n)_{n \in \ints}$ in $Y$. For all $n \in \ints$, using Hahn-Banach, choose $\phi_n \in X^*$ with $||\phi_n||=1$ such that $||x_n||_X = \phi_n(x_n)$.

Let $y \in Y$, and $\epsilon > 0$. Choose $n \in \ints$ such that $||y - x_n||_X < \epsilon / 2$.
\[ ||y||_X \leq ||x_n||_X + ||y - x_n||_X \leq \phi_n(x_n) +\epsilon / 2 \leq \phi_n(y) +||\phi_n||_{X^*} \ ||x_n - y||_X + \epsilon / 2,\]
so $\phi_n(y) \geq ||y||_X -\epsilon$.
\hfill $\square$
\end{quote}

\subsection{Lemma}
\label{lem:positivedecomposition}
Let $X$ be a Banach lattice and let $\phi \in X^*$. Then there exist $\phi^+, \phi^- \in X^*$ such that
\begin{itemize}
\item[(i)] $\phi = \phi^+ - \phi^-$;
\item[(ii)] $\phi^+(x)\geq 0$ for all $x \in X, x \succeq 0$;
\item[(iii)] $\phi^-(x) \geq 0$ for all $x \in X, x \succeq 0$.
\end{itemize}

\begin{quote}
\textsc{Proof:}
\end{quote}

\subsection{Lemma}
\label{lem:integralzero}
If $X$ is a Banach lattice, and $x \in L^1(M;X)$, with $x \succeq 0$ on $M$, then
\[ \int_M x \ d \mu = 0 \]
if and only if
\[ x = 0, \quad \mu \textrm{-almost everywhere}.\]

\begin{quote}
\textsc{Proof:}
Trivially, if $x = 0,$ $\mu$-almost everywhere on $M$, then $\int_M x \ d \mu = 0$.

Because $x \in L^1(M;X)$, there exists a closed separable subspace $Y \subset X$ with $x(\cdot) \in Y$, $\mu$-almost everywhere.
By Lemma~\ref{lem:normingsequence}, there exists a countable sequence $(\phi_n)_{n \in \ints}$ in $X^*$ that is norming for $Y$. By Lemma~\ref{lem:positivedecomposition}, every $\phi_n$ can be written as $\phi_n = \phi_n^+ - \phi_n^-$, with $\phi_n^+, \phi_n^-$ positive linear functionals. 

Because $\int_M x \ d \mu = 0$,
\[ \int_M \phi_n^{\pm} \circ x \ d \mu = \phi_n^\pm \left( \int_M x \ d \mu \right)  = 0,\]
and because $\phi_n^{\pm} \circ x \geq 0$, this shows that 
\[ \phi_n^{\pm} \circ  x  = 0, \quad \mu \textrm{-almost everywhere, for all} \ n \in \ints.\]
Because we only consider a countable number of $\phi_n$, this implies that
\[ \phi_n \circ x = 0, \quad \textrm{for all} \ n \in \ints, \ \mu\textrm{-almost everywhere}.\]
Therefore $x = 0, \mu$-almost everywhere.
\hfill $\square$
\end{quote}

\end{document}